\documentclass[a4paper,reqno]{amsart}

\usepackage[american]{babel}
\usepackage[T1]{fontenc}
\usepackage[utf8]{inputenc}
\usepackage{amsmath,amssymb,amsfonts,amsthm,amstext}
\usepackage{dsfont,mathtools,enumerate,mathrsfs}
\usepackage{color}
\usepackage{microtype}
\usepackage[numbers,sort&compress]{natbib}
\usepackage{hyperref}
\usepackage{ifthen}

\newcommand{\dx}{\mathrm{d}}
\newcommand{\apply}[3][]{\left<#2,#3\right>\ifthenelse{\equal{#1}{}}{}{_{#1}}}
\newcommand{\scalar}[3][]{\left(#2\mid#3\right)\ifthenelse{\equal{#1}{}}{}{_{#1}}}
\renewcommand{\phi}{\varphi}

\newcommand{\setone}{\mathds{1}}
\newcommand{\eps}{\varepsilon}
\renewcommand{\rho}{\varrho}
\newcommand{\wto}{\rightharpoonup}

\DeclareMathOperator{\sgn}{sgn}
\DeclareMathOperator{\divi}{div}
\DeclareMathOperator*{\esssup}{ess\,sup}
\DeclareMathOperator*{\essinf}{ess\,inf}

\DeclareMathOperator{\id}{id}

\theoremstyle{definition}
\newtheorem{theorem}{Theorem}[section]
\newtheorem{proposition}[theorem]{Proposition}
\newtheorem{corollary}[theorem]{Corollary}
\newtheorem{lemma}[theorem]{Lemma}
\newtheorem{remark}[theorem]{Remark}
\newtheorem{example}[theorem]{Example}
\newtheorem{definition}[theorem]{Definition}

\title[Quasilinear Robin Problems on Lipschitz Domains]{Quasilinear Elliptic and Parabolic Robin Problems on Lipschitz Domains}
\date{March 8, 2011}
\author{Robin Nittka}
\address{Robin Nittka\\University of Ulm\\Institute of Applied Analysis\\89069 Ulm\\Germany}
\email{robin.nittka@uni-ulm.de}

\numberwithin{equation}{section}

\begin{document}
\begin{abstract}
	We prove H\"older continuity up to the boundary for solutions of
	quasi-linear degenerate elliptic problems in divergence form, not necessarily of variational type,
	on Lipschitz domains with Neumann and Robin boundary conditions. This includes the
	$p$-Laplace operator for all $p \in (1,\infty)$, but also operators with unbounded coefficients.
	Based on the elliptic result we show that the corresponding parabolic problem is
	well-posed in the space $\mathrm{C}(\overline{\Omega})$ provided that the coefficients
	satisfy a mild monotonicity condition. More precisely, we show that the realization of the
	elliptic operator in $\mathrm{C}(\overline{\Omega})$ is m-accretive and densely defined.
	Thus it generates a non-linear strongly continuous contraction semigroup on $\mathrm{C}(\overline{\Omega})$.
\end{abstract}

\maketitle

\section{Introduction}

Given a bounded Lipschitz domain $\Omega$ in $\mathds{R}^N$, we show that all weak
solutions of certain degenerate quasi-linear elliptic problems are
H\"older-continuous up to the boundary of $\Omega$, which generalizes the results in~\cite{Nittka09}
to non-linear equations. More precisely, this is true for
equations of the form
\begin{equation}\label{eq:IntrRobinEq}
	\left\{ \begin{aligned}
		-\divi A(x,u,\nabla u) + B(x,u,\nabla u) + \omega u & = 0 && \text{on } \Omega \\
		A(x,u,\nabla u) \cdot \nu + h(x,u) & = 0 && \text{on } \partial\Omega
	\end{aligned} \right.
\end{equation}
where $A\colon \Omega \times \mathds{R} \times \mathds{R}^N \to \mathds{R}^N$,
$B\colon \Omega \times \mathds{R} \times \mathds{R}^N \to \mathds{R}$
and $h\colon \partial\Omega \times \mathds{R} \to \mathds{R}$ are measurable functions
such that there exist constants $1 < p < \infty$ and $0 < \nu \le \mu$
and non-negative functions $\psi_1$, $\psi_2$ and $\psi_3$ satisfying
\begin{equation}\label{eq:structure}
	\begin{aligned}
		z A(x,u,z) & \ge \nu |z|^p - \psi_1(x) |u|^p - \psi_1(x) \\
		|A(x,u,z)| & \le \mu |z|^{p-1} + \psi_2(x) |u|^{p-1} + \psi_2(x) \\
		|B(x,u,z)| & \le \psi_3(x) |z|^{p-1} + \psi_1(x) |u|^{p-1} + \psi_1(x) \\
		|h(x,u)| & \le \psi_4(x) |u|^{p-1} + \psi_4(x)
	\end{aligned}
\end{equation}
for all $x \in \Omega$, $u \in \mathds{R}$ and $z \in \mathds{R}^N$, and such that
\begin{align*}
	\psi_1 & \in L^{\frac{N}{p-\eps}}(\Omega), &
	\psi_2 & \in L^{\frac{N}{p-1}}(\Omega), &
	\psi_3 & \in L^{\frac{N}{1-\eps}}(\Omega), &
	\psi_4 & \in L^{\frac{N-1}{p-1-\eps}}(\partial\Omega)
	&& \text{if } p < N, \\
	\psi_1 & \in L^{\frac{N}{N-\eps}}(\Omega), &
	\psi_2 & \in L^{\frac{N}{N-1-\eps}}(\Omega), &
	\psi_3 & \in L^{\frac{N}{1-\eps}}(\Omega), &
	\psi_4 & \in L^{\frac{N-1}{N-1-\eps}}(\partial\Omega)
	&& \text{if } p = N, \\
	\psi_1 & \in L^1(\Omega), &
	\psi_2 & \in L^{\frac{p}{p-1}}(\Omega), &
	\psi_3 & \in L^p(\Omega), &
	\psi_4 & \in L^1(\partial\Omega)
	&& \text{if } p > N
\end{align*}
with some $\eps \in (0,1)$.
In particular, the $\Delta_p$-equation with Robin boundary conditions is included for all
$p \in (1,\infty)$, together with a large variety of lower order perturbations, including unbounded
coefficients. Also we do not require a variational structure of the equation. We refer
to~\cite{JLM01} and references therein for a short account on applications for the $\Delta_p$-operator.
The assumptions on the coefficients are optimal in that the regularity assumptions on $A$ and $B$ are the right
ones within the class of $L^p$-functions for results about interior regularity~\cite{Serrin}, see also~\cite{Lie91}.

The other main result of this article is that under suitable conditions on the coefficients the unique solution of
\begin{equation}\label{eq:IntrParRobin}
	\left\{ \begin{aligned}
		u_t(t,x) - \divi a(x,\nabla u(t,x)) + b(x,u(t,x)) & = 0 && t > 0, \; x \in \Omega \\
		a(x,\nabla u(t,x)) \cdot \nu + h(x,u(x)) & = 0 && t > 0, \; x \in \partial\Omega \\
		u(0,x) & = u_0(x) && x \in \Omega,
	\end{aligned} \right.
\end{equation}
which we define in an $L^2$-sense, is continuous on the parabolic cylinder
$[0,\infty) \times \overline{\Omega}$ whenever $u_0 \in \mathrm{C}(\overline{\Omega})$.
More precisely, we show that the corresponding elliptic operator is m-accretive on $\mathrm{C}(\overline{\Omega})$
and thus generates a non-linear contraction $\mathrm{C}_0$-semigroup on $\mathrm{C}(\overline{\Omega})$.
In order to obtain this result we have to assume that $a$, $b$ and $h$ are Carath\'eodory functions that such
$A(x,u,z) \coloneqq a(x,z)$, $B(x,u,z) \coloneqq b(x,u)$ and $h$
satisfy the above conditions~\eqref{eq:structure}. Moreover, we assume the following
mild monotonicity assumptions,
\begin{equation}\label{eq:IntrMonAss}
	\left\{ \begin{aligned}
		(z_1 - z_2) \, (a(x,z_1) - a(x,z_2)) & \ge 0, \\
		(u_1 - u_2) \, (b(x,u_1) - b(x,u_2)) & \ge 0, \\
		(u_1 - u_2) \, (h(x,u_1) - h(x,u_2)) & \ge 0,
	\end{aligned} \right.
\end{equation}
which are much weaker than the standard monotonicity assumptions as considered for example in~\cite{Min07}.
In particular, \eqref{eq:IntrMonAss} includes the $\Delta_p$-operator for every $p \in (1,\infty)$.
Thus for all $p \in (1,\infty)$ the problem
\[
	\left\{ \begin{aligned}
		u_t(t,x) - \Delta_p u(t,x) + b_0(x) |u|^{p-2} u & = 0 && t > 0, \; x \in \Omega \\
		|\nabla u(t,x)|^{p-2} \tfrac{\partial u(t,x)}{\partial \nu} + h_0(x) |u|^{p-2} u & = 0 && t > 0, \; x \in \partial\Omega \\
		u(0,x) & = u_0(x) && x \in \Omega
	\end{aligned} \right.
\]
is well-posed in $\mathrm{C}(\overline{\Omega})$ if $b_0 \in L^\infty(\Omega)$
and $h_0 \in L^\infty(\partial\Omega)$ are nonnegative. We also obtain a similar result
for Wentzell-Robin boundary conditions. All of these results are based on the
the author's PhD thesis~\cite{PhD}.

\smallskip

There are several good reasons to study elliptic and parabolic equations in $\mathrm{C}(\overline{\Omega})$.
On the one hand, for Dirichlet boundary conditions it is the natural space to formulate the boundary conditions, cf.~\cite{Are00,AC10,AS11},
and it is the natural space for maximum principles. But also for Neumann or Robin boundary conditions
the space $\mathrm{C}(\overline{\Omega})$ is nicer in some respects than the $L^q$-spaces
with $q \in [1,\infty)$. For instance, the composition operator $f \mapsto g \circ f$ is locally Lipschitz
continuous in $\mathrm{C}(\overline{\Omega})$ whenever $g$ is locally Lipschitz continuous, but in general fails
to map $L^q(\Omega)$ into $L^q(\Omega)$, so rapidly growing non-linear perturbations can more easily be handled
in $\mathrm{C}(\overline{\Omega})$ than in $L^q(\Omega)$.

The results of this article are new regarding several aspects. Our results are valid for bounded Lipschitz domains,
which form a strictly larger class than the strong Lipschitz domains, i.e., the domains that are locally
the epigraph of a Lipschitz continuous functions. For example, the physically relevant example of
the topologically regularized union of two crossing beams is a Lipschitz domain, but not a strong Lipschitz domain~\cite[\S 7.3]{HDR09}.
The class of Lipschitz domains has been studied a lot recently, see~\cite{HDR09,HDR10,GGKR02,KNR06}, to name only a few articles.

H\"older continuity of solutions of linear equations with Robin boundary conditions on Lipschitz domains has been extensively studied,
sometimes only in special cases, see for example~\cite{Warma06,Nittka09,GR01,Fuk67,BH91}.
The main elliptic result of this article, Theorem~\ref{thm:HoelderRobin},
seems to be new in the non-linear case even for smooth domains, but compare~\cite{GR06,Lie83}
for corresponding results under more restrictive assumptions on the coefficients,
which are obtained by different methods.
The linear parabolic problem has been studied in~\cite{Nittka09,Warma06} in terms of semigroups. The non-linear
case seems to be new, but see~\cite{GR10} for results under stronger regularity assumptions.
Our result seems to be particularly interesting because we neither assume that
the corresponding elliptic operator has a variational structure nor that it is strongly monotone.

\smallskip

The article is structured as follows. After introducing some notational conventions and basic properties
of Lipschitz domains in Section~\ref{sec:notation}, we show in Section~\ref{sec:neumann} that every
solution of~\eqref{eq:IntrRobinEq} for $h=0$ is H\"older continuous, thus proving the main elliptic
result for Neumann boundary conditions. The proof
is based on a reflection argument that the author has used already for the linear case~\cite{Nittka09}. The
general idea is much older, compare for example~\cite[Section~2.4.3]{Troi87}, but has apparently not
been exploited to this extent before.

In Section~\ref{sec:robin} we obtain a priori estimates for the Robin problem.
We use Moser's iteration in a similar manner as in~\cite{DD09}. Our result is more general than those in~\cite{DD09}
in that we allow general quasi-linear operators, but less general in that we restrict ourselves to Lipschitz domains.
Combining these a priori estimates with our main result for problems with Neumann boundary conditions
we extend the regularity result to general $h$.

Finally, in Section~\ref{sec:par} we make use of the elliptic theory in order to show that the parabolic problem with Robin boundary conditions
is well-posed in the space of continuous functions. The result is based on non-linear semigroup theory. Following the ideas in~\cite{AMPR03},
we are able to apply our methods also to equations with Wentzell-Robin boundary conditions.
We do not have to assume that the $L^2$-realization of the operator
is a subdifferential, i.e., we do not assume that the corresponding elliptic problem has a variational formulation.

\section{Notation and preliminaries}\label{sec:notation}

Throughout the article we follow the convention that constants denoted by $c$
are allowed to vary from one line to the next and may depend freely
on the parameters $N$, $\Omega$, $p$, $\eps$ and upper bounds for $\|\psi_1\|$, $\|\psi_2\|$, $\|\psi_3\|$ and $\|\psi_4\|$
in their respective spaces, as introduced in~\eqref{eq:structure}.
Any additional dependencies are explicitly indicated by subscripts.

\begin{definition}\label{def:Lipschitz}
	We say that an open set $\Omega \subset \mathds{R}^N$ is a \emph{Lipschitz domain} if
	for every $x \in \partial\Omega$ there exists an $\mathds{R}^N$-neighborhood $V$ of $x$
	and a bi-Lipschitz mapping $\psi$ from $V$ onto $(-1,1)^N$ such that
	$\psi(V \cap \Omega) = (-1,1)^{N-1} \times (0,1)$, i.e., $\psi$ is invertible and both
	$\psi$ and $\psi^{-1}$ are Lipschitz continuous.
	One says that $\overline{\Omega}$ is an $N$-dimensional, bounded Lipschitz submanifold
	of $\mathds{R}^N$ with boundary.
\end{definition}

\begin{remark}[{\cite[\S 1.2.1]{Grisvard}}]
	Every domain with Lipschitz boundary is a Lipschitz domain,
	but the converse fails.
\end{remark}

If $\Omega$ is a Lipschitz domain,
then it is an extension domain~\cite[Theorem~7.25]{GT} which implies
that $\mathrm{C}^\infty(\overline{\Omega})$ is dense in $W^{1,p}(\Omega)$
for $p \in [1,\infty)$ and that the usual Sobolev embeddings hold.
The natural measure on the boundary of a Lipschitz domain is the $(N-1)$-dimensional
Hausdorff measure. It is the unique measure for which the divergence theorem
holds. This fact is easily transported from the reference domain $(-1,1)^{N-1} \times (0,1)$
to $\Omega$, see also~\cite[\S 5.8]{EG}.
We agree that integrals over the boundary of a Lipschitz domain are always to be understood to be
taken with respect to the $(N-1)$-dimensional Hausdorff measure, which we denote by $\sigma$ if
the need arises.

Let $1 \le q \le \infty$.
If we write $\|u\|_{L^q(\Omega)}$, where $u$ is a measurable function on $\Omega$, we regard that
expression to equal infinity if $u \not\in L^q(\Omega)$. For convenience we use the notation
\[
	\|u\|_{L^q(\partial\Omega)} \coloneqq \bigl\| u|_{\partial\Omega} \bigr\|_{L^q(\partial\Omega)}
\]
for functions $u \in W^{1,p}(\Omega)$, which admit a trace $u|_{\partial\Omega} \in L^p(\partial\Omega)$,
and again we define this expression to equal infinity if $u|_{\partial\Omega} \not\in L^q(\partial\Omega)$.

We will need a change of variables formula for boundary integrals. In order to prove it,
we start with a few facts about the derivative of bi-Lipschitz mappings.

\begin{lemma}\label{lem:jacest}
	Let $U \subset \mathds{R}^n$ be open and let $\psi\colon U \to \mathds{R}^m$
	be a bi-Lipschitz mapping, $m \ge n$.
	Then $\psi$ is differentiable almost everywhere and the Jacobian
	$J\psi(x) \coloneqq (\det(\psi'(x)^T \psi'(x)))^{1/2}$ of $\psi$
	satisfies $\alpha \le J\psi \le \beta$ almost everywhere, where the constants
	$\alpha$ and $\beta$ depend only on $n$, $m$ and the Lipschitz constants of $\psi$ and
	$\psi^{-1}$. If $m=n$, then $\psi'$ is invertible almost everywhere with
	uniformly bounded inverse and $\alpha \le |\det \psi'| \le \beta$.
\end{lemma}
\begin{proof}
	Rademacher's theorem asserts that $\psi$ is differentiable almost everywhere.
	More precisely, the entries of $\psi'$ are essentially bounded by the Lipschitz
	constant of $\psi$, which proves the upper estimate for $J\psi$.
	As for the lower estimate, let $L > 0$ be the Lipschitz constant of $\psi^{-1}$, so that
	$|\psi(y) - \psi(x)| \ge L^{-1} |x-y|$ for all $x,y \in U$.
	If $x$ is a point of differentiability of $\psi$, then
	\[
		L^{-1} |tv| \le |\psi(x+tv) - \psi(x)| = |t \psi'(x)v + o(t)|.
	\]
	For $t \to 0$ we obtain that $|\psi'(x)v| \ge L^{-1} |v|$ almost everywhere for all $v \in \mathds{R}^n$,
	hence
	\[
		|v| \, |\psi'(x)^T \psi'(x)v|
			\ge \scalar{\psi'(x)^T \psi'(x)v}{v}
			= |\psi'(x) v|^2
			\ge L^{-2} |v|^2.
	\]
	Thus the eigenvalues of $\psi'(x)^T \psi'(x)$ can be bounded from below by $L^{-2}$, which shows that
	$J\psi \ge L^{-n}$ holds almost everywhere.

	Finally, if $m=n$, then the chain rule~\cite[Theorem~2.2.2]{Ziemer} implies that
	\[
		(\psi^{-1})'(\psi(x)) = (\psi'(x))^{-1}.
	\]
	Since the entries of $(\psi^{-1})'$
	are essentially bounded, this proves that $(\psi')^{-1}$ is uniformly bounded outside
	a set of measure zero. The estimate for the determinant follows from $J\psi = |\det \psi'|$.
\end{proof}

We can now prove the following change of variables formula for boundary integrals.

\begin{lemma}\label{lem:bdr_cov}
	Let $\Omega_1$ and $\Omega_2$ be Lipschitz domains and $\psi\colon \Omega_2 \to \Omega_1$
	a bi-Lipschitz function. Then $\psi$ has a unique extension to $\overline{\Omega}_1$,
	and $\psi(\partial\Omega_2) = \partial\Omega_1$, where we identify $\psi$ with its extension.
	In this situation, there exists a measurable function
	$m\colon \Omega_2 \to (0, \infty)$, which is unique up to nullsets,
	such that
	\[
		\int_{\partial\Omega_1} g = \int_{\partial\Omega_2} (g \circ \psi) m
	\]
	for all positive measurable functions $g$ on $\partial\Omega_1$ and hence for all
	integrable functions. Moreover, $0 < \alpha \le m \le \beta$ almost everywhere with
	constants $\alpha$ and $\beta$ that depend only on $\psi$, $\Omega_1$ and $\Omega_2$.
\end{lemma}

\begin{proof}
	The assertions about $\psi$ and the uniqueness of $m$ are clear.
	In order to show the existence of $m$, fix $y \in \partial\Omega_2$ and define $x \coloneqq \psi(y) \in \partial\Omega_1$.
	Fix neighborhoods $V_2$ of $y$ and $V_1$ of $x$ such that there exist bi-Lipschitz transformations
	$\psi_1\colon V_1 \to (-1,1)^N$ and $\psi_2\colon V_2 \to (-1,1)^N$ as in Definition~\ref{def:Lipschitz}.
	Without loss of generality we pick $V_2$ so small that $\psi(V_2) \subset V_1$.
	Write
	\[
		\phi_i \coloneqq \psi_i|_{\partial\Omega_i \cap V_i}\colon \partial\Omega_i \cap V_i \to H \coloneqq (-1,1)^{N-1} \times \{0\}.
	\]
	Then the bi-Lipschitz function $\phi \coloneqq \phi_1 \circ \psi \circ \phi_2^{-1}\colon H \to H$ is the local representation
	of $\psi|_{\partial\Omega_2}$, and we regard $H$ as an open subset of $\mathds{R}^{N-1}$.
	Then for every positive measurable function $f$ on $H$ we obtain that
	\[
		\int_H f = \int_H (f \circ \phi) |\det \phi'|
	\]
	by the change of variables formula for Lipschitz functions~\cite[\S 3.3.3]{EG}. Also by the change
	of variables formula we obtain that
	\[
		\int_H f \, J\phi_i^{-1} = \int_{\partial\Omega_i} f \circ \phi_i
	\]
	for $i=1,2$, where $J\phi_i^{-1}$ denotes the Jacobian of $\phi_i^{-1}$ as a mapping from $H$ into $\mathds{R}^N$.
	Combining these formulas, we see that
	\begin{align*}
		\int_{\partial\Omega_1} g
			& = \int_H (g \circ \phi_1^{-1}) \, J\phi_1^{-1}
			= \int_H (g \circ \phi_1^{-1} \circ \phi) (J\phi_1^{-1} \circ \phi) |\det\phi'| \\
			& = \int_{\partial\Omega_2} (g \circ \phi_1^{-1} \circ \phi \circ \phi_2) \frac{J\phi_1^{-1} \circ \phi \circ \phi_2}{J\phi_2^{-1} \circ \phi_2} |\det \phi' \circ \phi_2|
	\end{align*}
	for all positive measurable functions $g$ on $\partial\Omega_1$. By
	Lemma~\ref{lem:jacest}, this implies the claim.
\end{proof}

\section{Elliptic Neumann problems}\label{sec:neumann}

Let $\Omega \subset \mathds{R}^N$ be a Lipschitz domain and $p \in (1,\infty)$.
We prove that all weak solutions $u \in W^{1,p}(\Omega)$ of
\begin{equation}\label{eq:NeumannEq}
	\left\{ \begin{aligned}
		-\divi A(x,u,\nabla u) + B(x,u,\nabla u) & = f(x) - (\divi F)(x) && \text{on } \Omega \\
		A(x,u,\nabla u) \cdot \nu & = g(x) + (F \cdot \nu)(x) && \text{on } \partial\Omega
	\end{aligned} \right.
\end{equation}
are H\"older continuous, provided $A$ satisfies the structure conditions~\eqref{eq:structure} and
$f$, $F$ and $g$ are in appropriate Lebesgue spaces, namely
\begin{equation}\label{eq:fFass}
	\left\{ \begin{aligned}
		f & \in L^{\frac{N}{p-\eps}}(\Omega), &
		F & \in L^{\frac{N}{p-1}}(\Omega;\mathds{R}^N), &
		g & \in L^{\frac{N-1}{p-1}}(\partial\Omega) &
		&& \text{if } p < N, \\
		f & \in L^{\frac{N}{N-\eps}}(\Omega), &
		F & \in L^{\frac{N}{N-1-\eps}}(\Omega;\mathds{R}^N), &
		g & \in L^{\frac{N-1}{N-1-\eps}}(\partial\Omega) &
		&& \text{if } p = N, \\
		f & \in L^1(\Omega), &
		F & \in L^{\frac{p}{p-1}}(\Omega;\mathds{R}^N), &
		g & \in L^1(\partial\Omega) &
		&& \text{if } p > N.
	\end{aligned} \right.
\end{equation}
We could absorb $f$ and $F$ (but not $g$) into the coefficients, but for the application
we have in mind it turns out to be more convenient to write them down explicitly.
By convention, $\|f\|$, $\|F\|$ and $\|g\|$ will always refer to the norms of $f$, $F$
and $g$ in the spaces indicated in~\eqref{eq:fFass}.

\begin{definition}
	We say that $u \in W^{1,p}(\Omega)$ is a weak solution of~\eqref{eq:NeumannEq} if
	\begin{equation}\label{eq:weakform}
		\int_\Omega \nabla\eta \, A(x,u,\nabla u) + \int_\Omega \eta \, B(x,u,\nabla u)
			 = \int_\Omega \eta \, f + \int_\Omega \nabla\eta \, F + \int_{\partial\Omega} \eta \, g
	\end{equation}
	for all $\eta \in \mathrm{C}^\infty(\overline{\Omega})$. If~\eqref{eq:weakform} holds merely for all
	$\eta \in \mathrm{C}^\infty_c(\Omega)$, we say that $u$ is a solution of the equation given by
	the first line of~\eqref{eq:NeumannEq}, without any boundary conditions.
\end{definition}

\begin{remark}
	A function $u \in W^{1,p}(\Omega)$ is a weak solution of~\eqref{eq:NeumannEq} if and only if~\eqref{eq:weakform}
	holds for all $\eta \in W^{1,p}(\Omega)$ since $\mathrm{C}^\infty(\overline{\Omega})$ is dense in $W^{1,p}(\Omega)$
	and all expressions in~\eqref{eq:weakform} are continuous as $\eta$ varies in $W^{1,p}(\Omega)$, compare
	Proposition~\ref{prop:monop}, where an even stronger assertion is proved. In what follows, we will use this fact frequently.
\end{remark}

We deduce boundary regularity from the following interior regularity result, which is
an immediate consequence of results due to Serrin.
\begin{theorem}[{\cite[\S 1.1, \S 1.4, \S 1.5]{Serrin}}]\label{thm:serrin}
	Let $A$ and $B$ satisfy the structure conditions~\eqref{eq:structure} and
	let $f$ and $F$ be as in~\eqref{eq:fFass}.
	Then there exists $\alpha \in (0,1)$
	such that every weak solution $u \in W^{1,p}(\Omega)$ of
	\begin{equation}\label{eq:interioreq}
		-\divi A(x,u,\nabla u) + B(x,u,\nabla u) = f(x) - (\divi F)(x) \text{ on } \Omega
	\end{equation}
	is in $\mathrm{C}^{0,\alpha}_{\mathrm{loc}}(\Omega)$. Moreover, for every
	relatively compact subdomain $\omega \subset \Omega$ there exists $c_{\alpha,\omega} \ge 0$ such that
	\begin{equation}\label{eq:serrin:est}
		\|u\|_{\mathrm{C}^{0,\alpha}(\omega)} \le c_{\alpha,\omega} \bigl( \|f\|^{\frac{1}{p-1}} + \|F\|^{\frac{1}{p-1}} + \|u\|_{L^p(\Omega)} \bigr) + c_{\alpha,\omega}
	\end{equation}
	holds for all weak solutions $u \in W^{1,p}(\Omega)$ of~\eqref{eq:interioreq}.
\end{theorem}

In order to apply Theorem~\ref{thm:serrin}, we extend the solutions of~\eqref{eq:NeumannEq}
locally along the boundary of $\Omega$ and show that the extension satisfies an elliptic equation
on the larger domain. Then interior regularity of the extended function implies boundary
regularity of the original function.

It is convenient to carry over the function to the reference domain $(-1,1)^{N-1} \times (0,1)$ and
to extend the resulting function on that domain.
As a first step, we show that the structural properties of the equation are preserved under bi-Lipschitz transformations.
\begin{proposition}\label{prop:transform}
	Let $\Omega_1$ and $\Omega_2$ be Lipschitz domains in $\mathds{R}^N$
	and let $\psi\colon \Omega_2 \to \Omega_1$ be a bi-Lipschitz bijection.
	Let $u \in W^{1,p}(\Omega_1)$ and $v \coloneqq u \circ \psi$.
	Given functions $A$, $B$, $f$, $F$ and $g$ as in~\eqref{eq:structure} and~\eqref{eq:fFass},
	define
	\begin{align*}
		\hat{A}(x,u,z) & \coloneqq (\psi'(x))^{-1} \, A(\psi(x), u, z \, \psi'(x)^{-1}) \; |\det\psi'(x)| \\
		\hat{B}(x,u,z) & \coloneqq B(\psi(x), u, z \, \psi'(x)^{-1}) \; |\det\psi'(x)|
	\end{align*}
	for $x \in \Omega_2$, $u \in \mathds{R}$ and row vectors $z \in \mathds{R}^N$.
	Moreover, let
	$\hat{f} \coloneqq (f \circ \psi) \, |\det\psi'|$, $\hat{F} \coloneqq (\psi')^{-1} \, (F \circ \psi) \, |\det\psi'|$
	and $\hat{g} \coloneqq (g \circ \psi) m$ with $m$ as in Lemma~\ref{lem:bdr_cov}.
	\begin{enumerate}[(a)]
	\item\label{eq:transform:reg}
		The function $v$ is in $W^{1,p}(\Omega_2)$ with $\nabla v = (\nabla u \circ \psi) \, \psi'$
		almost everywhere, and the functions $\hat{f}$,
		$\hat{F}$ and $\hat{g}$ are in Lebesgue spaces with the same exponent as $f$,
		$F$ and $g$, respectively. More precisely,
		$\|\hat{f}\| \le c_\psi \|f\|$, $\|\hat{F}\| \le c_\psi \|F\|$
		and $\|\hat{g}\| \le c_\psi \|g\|$.
	\item\label{eq:transform:struct}
		The functions $\hat{A}$ and $\hat{B}$ satisfy the structure conditions~\eqref{eq:structure} on $\Omega_2$,
		where the parameters depend only on $\psi$ and the parameters for $A$ and $B$.
	\item\label{eq:transform:neweq}
		If $u \in W^{1,p}(\Omega_1)$ satisfies
		\begin{equation}\label{eq:transform:eq1}
			\int_{\Omega_1} \nabla\eta \, A(x,u,\nabla u) + \int_{\Omega_1} \eta \, B(x,u,\nabla u)
				= \int_{\Omega_1} \eta \, f + \int_{\Omega_1} \nabla \eta \, F + \int_{\partial\Omega_1} \eta \, g
		\end{equation}
		for all $\eta \in \mathrm{C}^\infty_c(\Omega_1 \cup \Gamma)$
		with some relatively open set $\Gamma \subset \partial\Omega_1$ then $v \in W^{1,p}(\Omega_2)$ satisfies
		\begin{equation}\label{eq:transform:eq2}
			\int_{\Omega_2} \nabla\tilde{\eta} \, \hat{A}(x,v,\nabla v) + \int_{\Omega_2} \tilde{\eta} \, \hat{B}(x,v,\nabla v)
				= \int_{\Omega_2} \tilde{\eta} \, \hat{f} + \int_{\Omega_2} \nabla \tilde{\eta} \, \hat{F} + \int_{\partial\Omega_2} \tilde{\eta} \, \hat{g}
		\end{equation}
		for all $\tilde{\eta} \in \mathrm{C}^\infty_c(\Omega_2 \cup \psi^{-1}(\Gamma))$.
	\end{enumerate}
\end{proposition}
\begin{proof}\allowdisplaybreaks
	The assertions in~\eqref{eq:transform:reg} follow from
	the chain rule for Sobolev functions~\cite[Theorem~2.2.2]{Ziemer},
	the change of coordinates formula for Lipschitz transformations~\cite[\S 3.3.3]{EG}
	and Lemma~\ref{lem:bdr_cov}. We have also used that $|\det\psi'|$ is bounded from below
	and above, see Lemma~\ref{lem:jacest}.

	In order to check~\eqref{eq:transform:struct} we fix $u \in \mathds{R}$
	and a row vector $z \in \mathds{R}^N$. Then
	\begin{align*}
		z \, \hat{A}(x,u,z)
			& = z (\psi'(x))^{-1} \, A(\psi(x), u, z \, (\psi'(x))^{-1}) \; |\det\psi'(x)| \\
			& \ge \Bigl( \nu \bigl| z \; \psi'(x)^{-1} \bigr|^p - \psi_1(\psi(x)) \, |u|^p - \psi_1(\psi(x)) \Bigr)  \; |\det\psi'(x)| \\
			& \ge \essinf_{\Omega_2} |\det\psi'| \cdot \Bigl( \frac{\nu}{\esssup_{\Omega_2}\|\psi'\|^p} |z|^p - \hat{\psi}_1 \, |u|^p - \hat{\psi}_1 \Bigr)
	\end{align*}
	with the function $\hat{\psi}_1 \coloneqq \psi_1 \circ \psi$ possessing
	the same degree of integrability as $\psi_1$. Similarly,
	\begin{align*}
		|\hat{A}(x,u,z)| & \le \esssup_{\Omega_2} |\det\psi'| \; s \bigl( \mu s^{p-1} \, |z|^{p-1} + \hat{\psi}_2 |u|^{p-1} + \hat{\psi}_2 \bigr) \\
		|\hat{B}(x,u,z)| & \le \esssup_{\Omega_2} |\det\psi'| \; \bigl( s^{p-1} \hat{\psi}_3 |z|^{p-1} + \hat{\psi}_1 |u|^{p-1} + \hat{\psi}_1 \bigr)
	\end{align*}
	where $s \coloneqq \esssup_{\Omega_2} \|(\psi')^{-1}\|$, $\hat{\psi}_2 \coloneqq \psi_2 \circ \psi$
	and $\hat{\psi}_3 \coloneqq \psi_3 \circ \psi$. Hence $\hat{A}$ and $\hat{B}$ satisfy~\eqref{eq:structure} on $\Omega_2$.

	As for~\eqref{eq:transform:neweq}, let $\Gamma \subset \partial\Omega_1$ be relatively open
	and assume that $u$ satisfies~\eqref{eq:transform:eq1} for all $\eta \in \mathrm{C}^\infty_c(\Omega_1 \cup \Gamma)$.
	Then by denseness~\eqref{eq:transform:eq1} is fulfilled for all $\eta \in W^{1,p}_0(\Omega_1 \cup \Gamma)$,
	the closure of $\mathrm{C}^\infty_c(\Omega_1 \cup \Gamma)$ in $W^{1,p}(\Omega)$.
	Let $\tilde{\eta}$ be in $\mathrm{C}^\infty_c(\Omega_2 \cup \psi^{-1}(\Gamma))$ and write $\eta \coloneqq \tilde{\eta} \circ \psi^{-1}$.
	Then by a standard smoothing argument we obtain that $\eta \in W^{1,p}_0(\Omega_1 \cup \Gamma)$. Moreover,
	\[
		\nabla\eta = (\nabla\tilde{\eta} \circ \psi^{-1}) (\psi^{-1})' = (\nabla\tilde{\eta} \circ \psi^{-1}) (\psi' \circ \psi^{-1})^{-1}
	\]
	by the chain rule. Hence
	\begin{align*}
		& \int_{\Omega_2} \nabla\tilde{\eta}(x) \, \hat{A}(x,v(x),\nabla v(x))
				+ \int_{\Omega_2} \tilde{\eta}(x) \, \hat{B}(x,v(x),\nabla v(x)) \\
			& \qquad = \int_{\Omega_2} \nabla\eta(\psi(x)) \, A(\psi(x), u(\psi(x)), \nabla u(\psi(x))) \, |\det\psi'(x)| \\
				& \qquad\qquad  + \int_{\Omega_2} \eta(\psi(x)) \, B(\psi(x), u(\psi(x)), \nabla u(\psi(x))) \, |\det\psi'(x)| \\
			& \qquad = \int_{\Omega_1} \nabla\eta(x) \, A(x,u(x),\nabla u(x))
				+ \int_{\Omega_1} \eta(x) \, B(x,u(x), \nabla u(x)) \\
			& \qquad = \int_{\Omega_1} \eta \, f
				+ \int_{\Omega_1} \nabla\eta \, F
				+ \int_{\partial\Omega_1} \eta \, g \\
			& \qquad = \int_{\Omega_2} \tilde{\eta} \, (f \circ \psi) \, |\det\psi'|
				+ \int_{\Omega_2} \nabla\tilde{\eta} \, (\psi')^{-1} \, (F \circ \psi) \, |\det\psi'|
				+ \int_{\partial\Omega_2} \tilde{\eta} \, (g \circ \psi) m \\
			& \qquad = \int_{\Omega_2} \tilde{\eta} \, \hat{f} + \int_{\Omega_2} \nabla\tilde{\eta} \, \hat{F} + \int_{\partial\Omega_2} \tilde{\eta} \, \hat{g},
		\end{align*}
	where by the change of variables formula for integrals over the domain as well as over its boundary, see Lemma~\ref{lem:bdr_cov}.
\end{proof}

\begin{remark}\label{rem:transform}
	In part~\eqref{eq:transform:neweq} of Proposition~\ref{prop:transform}, the values of $g$ on $\partial\Omega_1 \setminus \Gamma$
	do not appear in the assertions. Hence we can apply the result also if $g$ is given only on $\Gamma$ and is left undefined on the
	remaining part of $\partial\Omega_1$ by artificially defining $g \coloneqq 0$ on $\partial\Omega_1 \setminus \Gamma$.
\end{remark}

We will also use the following representation of boundary integrals as integrals over the domain itself.
\begin{lemma}\label{lem:bdrtoint}
	Let $\Omega' \coloneqq (-1,1)^N$ and $H \coloneqq (-1,1)^N \times \{0\}$,
	and let $g \in L^q(H)$, $q \in [1,\infty)$. There exist functions $k$ and $K$ in the spaces
	\[
		\left\{ \begin{aligned}
			k & \in L^{\frac{qN}{N-1}}(\Omega'), &
			K & \in L^{\frac{qN}{N-1}}(\Omega';\mathds{R}^N) && \text{if } q > 1 \\
			k & \in L^{\frac{N}{N-1+\eps}}(\Omega'), &
			K & \in L^{\frac{N}{N-1+\eps}}(\Omega';\mathds{R}^N) && \text{if } q = 1,
		\end{aligned} \right.
	\]
	where $\eps \in (0,1)$ is arbitrary, such that $k$ and $K$ satisfy
	\[
		\int_H \eta g = \int_{\Omega'} \eta k + \int_{\Omega'} \nabla \eta \, K
	\]
	for all $\eta \in \mathrm{C}^\infty_c(\Omega')$. Moreover,
	$\|k\| + \|K\| \le c_{q,\eps} \|g\|_{L^q(H)}$.
\end{lemma}
\begin{proof}
	Define the linear functional $\phi\colon \mathrm{C}^\infty_c(\Omega') \to \mathds{R}$
	by $\phi(\eta) \coloneqq \int_H \eta g$. If $q > 1$, then
	\[
		|\phi(\eta)| \le \|\eta\|_{L^{\frac{q}{q-1}}(H)} \|g\|_{L^q(H)}
			\le c_q \|\eta\|_{W^{1,\frac{qN}{qN-N+1}}(\Omega')} \|g\|_{L^q(H)}
	\]
	by the Sobolev embedding theorems. Similarly, if $q=1$, then
	\[
		|\phi(\eta)| \le \|\eta\|_{L^\infty(H)} \|g\|_{L^1(H)}
			\le c_\eps \|\eta\|_{W^{1,\frac{N}{1-\eps}}(\Omega')} \|g\|_{L^1(H)}
	\]
	for every $\eps > 0$. Hence $\phi$ extends to a continuous
	linear functional on $W^{1,\frac{qN}{qN-N+1}}_0(\Omega')$ or $W^{1,\frac{N}{1-\eps}}_0(\Omega')$,
	respectively, which implies the claim, see~\cite[\S 4.3]{Ziemer}.
\end{proof}

We now prove the main result of this section: every weak solution of~\eqref{eq:NeumannEq} is H\"older continuous
up to the boundary of $\Omega$.
\begin{theorem}\label{thm:HoelderNeumann}
	Let $\Omega$ be a Lipschitz domain and assume~\eqref{eq:structure} and~\eqref{eq:fFass}.
	Then there exist $\alpha \in (0,1)$ and $c_\alpha \ge 0$ such that every weak solution $u \in W^{1,p}(\Omega)$
	of~\eqref{eq:NeumannEq} is in $\mathrm{C}^{0,\alpha}(\Omega)$ and satisfies
	\begin{equation}\label{eq:HoelderNeumann:est}
		\|u\|_{\mathrm{C}^{0,\alpha}(\Omega)}
			\le c_\alpha \bigl( \|f\|^{\frac{1}{p-1}} + \|F\|^{\frac{1}{p-1}} + \|g\|^{\frac{1}{p-1}} + \|u\|_{L^p(\Omega)} \bigr) + c_\alpha.
	\end{equation}
\end{theorem}
\begin{proof}
	Let $u \in W^{1,p}(\Omega)$ be a weak solution of~\eqref{eq:NeumannEq}.
	Let $z$ be in $\partial\Omega$ and fix $V$ and $\psi\colon V \to (-1,1)^N$
	as in Definition~\ref{def:Lipschitz}. We will show that $u$ is H\"older continuous in a
	neighborhood $V_z \subset V$ of $z$. Since $u$ solves~\eqref{eq:NeumannEq} we know in particular that
	\[
		\int_\Omega \nabla\eta \, A(x,u,\nabla u) + \int_\Omega \eta \, B(x,u,\nabla u)
			= \int_\Omega \eta \, f + \int_\Omega \nabla \eta \, F + \int_{\partial\Omega} \eta \, g
	\]
	for all $\eta \in \mathrm{C}^\infty_c(\Omega \cup (V \cap \partial\Omega))$.
	Write $\Omega_1 \coloneqq \psi(\Omega \cap V) = (-1,1)^{N-1} \times (0,1)$.
	Then by Proposition~\ref{prop:transform}, see also Remark~\ref{rem:transform},
	the function $v_1 \coloneqq u \circ \psi^{-1}$ satisfies
	\[
		\int_{\Omega_1} \nabla\eta \, \hat{A}_1(x,v_1,\nabla v_1) + \int_{\Omega_1} \eta \, \hat{B}_1(x,v_1,\nabla v_1)
			= \int_{\Omega_1} \eta \, \hat{f}_1 + \int_{\Omega_1} \nabla \eta \, \hat{F}_1 + \int_{\partial\Omega_1} \eta \, \hat{g}_1
	\]
	for all $\eta \in \mathrm{C}^\infty_c( (-1,1)^{N-1} \times [0,1) )$
	with functions $\hat{A}_1$, $\hat{B}_1$, $\hat{f}_1$, $\hat{F}_1$ and $\hat{g}_1$
	that satisfy the conditions~\eqref{eq:structure} and~\eqref{eq:fFass},

	Define the reflection $\psi_0\colon \mathds{R}^N \to \mathds{R}^N$
	by
	\[
		\psi_0(x_1,\dots,x_{N-1},x_N) \coloneqq (x_1,\dots,x_{N-1},-x_N).
	\]
	Then $v_2 \coloneqq v_1 \circ \psi_0 \in W^{1,p}(\Omega_2)$ satisfies
	\[
		\int_{\Omega_2} \nabla\eta \, \hat{A}_2(x,v_2,\nabla v_2) + \int_{\Omega_2} \eta \, \hat{B}_2(x,v_2,\nabla v_2)
			= \int_{\Omega_2} \eta \, \hat{f}_2 + \int_{\Omega_2} \nabla \eta \, \hat{F}_2 + \int_{\partial\Omega_2} \eta \, \hat{g}_2
	\]
	for all $\eta \in \mathrm{C}^\infty_c( (-1,1)^{N-1} \times (-1,0] )$ by Proposition~\ref{prop:transform}
	with functions $\hat{A}_2$, $\hat{B}_2$, $\hat{f}_2$, $\hat{F}_2$ and $\hat{g}_2$ that satisfy the conditions~\eqref{eq:structure} and~\eqref{eq:fFass}.

	Define $\Omega_0 \coloneqq (-1,1)^N$ and regard
	\begin{align*}
		v & \coloneqq v_1 \setone_{\Omega_1} + v_2 \setone_{\Omega_2}, \\
		\hat{A} & \coloneqq \hat{A}_1 \setone_{\Omega_1} + \hat{A}_2 \setone_{\Omega_2}, &
		\hat{B} & \coloneqq \hat{B}_1 \setone_{\Omega_1} + \hat{B}_2 \setone_{\Omega_2}, \\
		\hat{f} & \coloneqq \hat{f}_1 \setone_{\Omega_1} + \hat{f}_2 \setone_{\Omega_2}, &
		\hat{F} & \coloneqq \hat{F}_1 \setone_{\Omega_1} + \hat{F}_2 \setone_{\Omega_2}
	\end{align*}
	as functions on $\Omega_0$. Then $\hat{A}$, $\hat{B}$, $\hat{f}$, $\hat{F}$ and $\hat{g}$ satisfy the
	conditions~\eqref{eq:structure} and~\eqref{eq:fFass}. Moreover, $v \in W^{1,p}(\Omega_0)$
	and $\nabla v = \nabla v_1 \setone_{\Omega_1} + \nabla v_2 \setone_{\Omega_2}$.
	In fact, the Gauss-Green theorem~\cite[\S 5.8]{EG} shows that for all $\eta \in \mathrm{C}^\infty_c(\Omega_0)$ we have
	\begin{align*}
		\int_{\Omega_0} \nabla\eta \, v
			& = \int_{\Omega_1} \nabla\eta \, v_1 + \int_{\Omega_2} \nabla\eta \, v_2 \\
			& = \int_{\partial\Omega_1} \eta \, v_1 \, \nu_{\Omega_1} - \int_{\Omega_1} \eta \, \nabla v_1
				+ \int_{\partial\Omega_2} \eta \, v_2 \, \nu_{\Omega_2} - \int_{\Omega_2} \eta \, \nabla v_2 \\
			& = -\int_{\Omega_1} \eta \, \nabla v_1 - \int_{\Omega_2} \eta \, \nabla v_2
	\end{align*}
	since $v_1 = v_2$ on the intersection of the boundaries
	and the outer normals equal $\nu_{\Omega_1} = -e_N$ and $\nu_{\Omega_2} = e_N$ on that set,
	where $e_N$ denotes the $N$\textsuperscript{th} unit vector in $\mathds{R}^N$.
	In addition, $\eta$ vanishes on the remaining parts of $\partial\Omega_1$ and $\partial\Omega_2$ by assumption.

	Using in addition Lemma~\ref{lem:bdrtoint} we thus obtain that
	\begin{align*}
		& \int_{\Omega_0} \nabla\eta \, \hat{A}(x,v,\nabla v) + \int_{\Omega_0} \eta \, \hat{B}(x,v,\nabla v) \\
			& \quad = \int_{\Omega_1} \eta \, \hat{f}_1 + \int_{\Omega_1} \nabla\eta \, \hat{F}_1 + \int_{\partial\Omega_1} \eta \, \hat{g}_1
				+ \int_{\Omega_2} \eta \, \hat{f}_2 + \int_{\Omega_2} \nabla\eta \, \hat{F}_2 + \int_{\partial\Omega_2} \eta \, \hat{g}_2 \\
			& \quad = \int_{\Omega_0} \eta \, (\hat{f} + k) + \int_{\Omega_0} \nabla\eta \, (\hat{F} + K),
	\end{align*}
	for all $\eta \in \mathrm{C}^\infty_c(\Omega_0)$ with functions $k \in L^q(\Omega_0)$ and $K \in L^q(\Omega_0;\mathds{R}^N)$, where
	$q \coloneqq \frac{N}{p-1}$ for $p < N$, $q \coloneqq \frac{N}{N-1-\eps}$ for $p = N$ and $q \coloneqq \frac{p}{p-1}$ for $p > N$.
	Thus $v \in W^{1,p}(\Omega_0)$ is a weak solution of
	\[
		-\divi\hat{A}(x,v,\nabla v) + \hat{B}(x,v,\nabla v) = (\hat{f} + k) - \divi(\hat{F} + K) \text{ on } \Omega_0,
	\]
	where the coefficients $\hat{A}$ and $\hat{B}$ and the right hand side $\hat{f}+k$ and $\hat{F}+K$ satisfy the assumptions
	of Theorem~\ref{thm:serrin}. Consequently, there exists $\alpha \in (0,1)$ such that
	\begin{align*}
		\|v_1\|_{\mathrm{C}^{0,\alpha}( (-\frac{1}{2},\frac{1}{2})^{N-1} \times (0,\frac{1}{2}) )}
			& \le \|v\|_{\mathrm{C}^{0,\alpha}( (-\frac{1}{2},\frac{1}{2})^N )} \\
			& \le c_\alpha \bigl( \|\hat{f} + k\|^{\frac{1}{p-1}} + \|F+K\|^{\frac{1}{p-1}} + \|v\|_{L^p(\Omega_0)} \bigr) + c_\alpha \\
			& \le c_\alpha \bigl( \|f\|^\frac{1}{p-1} + \|F\|^{\frac{1}{p-1}} + \|g\|^{\frac{1}{p-1}} + \|u\|_{L^p(\Omega)} \bigr) + c_\alpha,
	\end{align*}
	see also part~\eqref{eq:transform:reg} of Proposition~\ref{prop:transform}. Since $u = v_1 \circ \psi$ and $\psi$ is Lipschitz
	continuous on $V$, we have shown that there exists a neighborhood $V_z \subset V$ of $z$ such that $u|_{V_z} \in \mathrm{C}^{0,\alpha}(V_z)$
	and
	\[
		\|u\|_{\mathrm{C}^{0,\alpha}(V_z)}
			\le c_\alpha \bigl( \|f\|^\frac{1}{p-1} + \|F\|^{\frac{1}{p-1}} + \|g\|^{\frac{1}{p-1}} + \|u\|_{L^p(\Omega)} \bigr) + c_\alpha.
	\]

	Since $\partial\Omega$ is compact, there exist finitely many $z_i \in \partial\Omega$ such that $\partial\Omega \subset \bigcup_{i=1}^m V_{z_i}$.
	Set $\omega_i \coloneqq V_{z_i} \cap \Omega$ and $\omega_0 \coloneqq \Omega \setminus \bigcup_{i=1}^m \omega_i \Subset \Omega$.
	Then $u|_{\omega_0} \in \mathrm{C}^{0,\alpha_0}(\omega_0)$ by Theorem~\ref{thm:serrin} and $u|_{\omega_i} \in \mathrm{C}^{0,\alpha_i}(\omega_i)$
	for $i=1,\dots,m$ by what we have just shown.
	Thus we have proved that $u \in \mathrm{C}^{0,\alpha}(\Omega)$ for $\alpha \coloneqq \min_{i=1,\dots,N} \alpha_i$,
	and more precisely we have shown that~\eqref{eq:HoelderNeumann:est} holds.
\end{proof}

\section{Elliptic Robin problems}\label{sec:robin}

Let $\Omega \subset \mathds{R}^N$ be a Lipschitz domain.
In this section we prove that all weak solutions $u \in W^{1,p}(\Omega)$ of
\begin{equation}\label{eq:RobinEq}
	\left\{ \begin{aligned}
		-\divi A(x,u,\nabla u) + B(x,u,\nabla u) + \omega u & = f(x) - (\divi F)(x) && \text{on } \Omega \\
		A(x,u,\nabla u) \cdot \nu + h(x,u) & = g(x) + (F \cdot \nu)(x) && \text{on } \partial\Omega
	\end{aligned} \right.
\end{equation}
are H\"older continuous, provided
$A$, $B$ and $h$ satisfy the structure conditions~\eqref{eq:structure},
$f$, $F$ and $g$ are as in~\eqref{eq:fFass} and $\omega$ is a nonnegative constant.
Like for Neumann boundary conditions, we say that $u \in W^{1,p}(\Omega)$ is a weak solution of~\eqref{eq:RobinEq} if
\[
	\int_\Omega \nabla\eta \, A(x,u,\nabla u) + \int_\Omega \eta \, B(x,u,\nabla u) + \int_{\partial\Omega} \eta \, h(x,u)
		= \int_\Omega \eta \, f + \int_\Omega \nabla\eta \, F + \int_{\partial\Omega} \eta \, g
\]
holds for all $\eta \in \mathrm{C}^\infty(\overline{\Omega})$ or, equivalently, for all $\eta \in W^{1,p}(\Omega) \cap L^2(\Omega)$.

We start by proving $L^q$-bounds for solutions of the slightly simpler equation
\begin{equation}\label{eq:RobinEqHom}
	\left\{ \begin{aligned}
		-\divi A(x,u,\nabla u) + B(x,u,\nabla u) + \omega u & = 0 && \text{on } \Omega \\
		A(x,u,\nabla u) \cdot \nu + h(x,u) & = 0 && \text{on } \partial\Omega
	\end{aligned} \right.
\end{equation}
via Moser's iteration technique. If we assume some strong monotonicity of the
coefficients, the a priori estimates could be obtained via an elegant
interpolation argument, see~\cite{PhD}. But in order to cover the general case
we have to use the iteration procedure instead.

We need the following lemma, whose easy proof we omit. But compare~\cite[Proposition~II.5.2]{Showalter}
for a similar argument.
\begin{lemma}\label{lem:eberlein}
	Let $X$, $Y$ and $Z$ be Banach spaces. Assume that $X$ is reflexive.
	Let $T\colon X \to Y$ be a compact linear operator and let $S\colon X \to Z$
	be an injective bounded linear operator.
	Then for every $\delta > 0$ there exists $c_\delta > 0$ such that
	\[
		\|Tx\|_Y \le \delta \|x\|_X + c_\delta \|Sx\|_Z
	\]
	holds for all $x \in X$.
\end{lemma}

\begin{proposition}\label{prop:Lqest}
	If $u \in W^{1,p}(\Omega)$ solves~\eqref{eq:RobinEqHom}, then
	\begin{equation}\label{eq:W1pest}
		\|u\|_{W^{1,p}(\Omega)} \le c \|u\|_{L^p(\Omega)} + c
	\end{equation}
	and
	\begin{equation}\label{eq:Lqest}
		\|u\|_{L^q(\Omega)} + \|u\|_{L^q(\partial\Omega)} \le c_q \|u\|_{L^p(\Omega)} + c_q
	\end{equation}
	for every $q \in [p,\infty)$.
\end{proposition}
\begin{proof}
	Let $u \in W^{1,p}(\Omega)$ be a weak solution of~\eqref{eq:RobinEqHom},
	and let $q \in [p,\infty)$ be arbitrary.
	Fix $\alpha \ge 1$ and define
	\[
		v_\alpha \coloneqq
				\bigl( (|u|+1)^{q-p+1} - 1 \bigr) \sgn(u) \setone_{\{|u| \le \alpha\}}
				+ \frac{(\alpha+1)^{q-p+1} - 1}{\alpha + 1} (|u|+1) \sgn(u) \setone_{\{|u| > \alpha\}}
	\]
	and
	\[
		w_\alpha \coloneqq
				(|u|+1)^{\frac{q}{p}} \setone_{\{|u| \le \alpha\}}
				+ (\alpha+1)^{\frac{q-p}{p}} (|u|+1) \setone_{\{|u| > \alpha\}}.
	\]
	Then by the chain rule~\cite[Theorem~2.1.11]{Ziemer} the functions $v_\alpha$ and $w_\alpha$ are in $W^{1,p}(\Omega)$ with weak derivatives
	\begin{equation}\label{eq:nablavalpha}
		\nabla v_\alpha = (q-p+1) (|u| + 1)^{q-p} \nabla u \, \setone_{\{|u| \le \alpha\}}
			+ \frac{(\alpha+1)^{q-p+1} - 1}{\alpha + 1} \nabla u \, \setone_{\{|u| > \alpha\}}
	\end{equation}
	and
	\[
		|\nabla w_\alpha| = \frac{q}{p} (|u| + 1)^{\frac{q-p}{p}} |\nabla u| \, \setone_{\{|u| \le \alpha\}}
			+ (\alpha + 1)^{\frac{q-p}{p}} |\nabla u| \, \setone_{\{|u| > \alpha\}},
	\]
	hence
	\begin{equation}\label{eq:nablawalpha}
		\Bigl( \frac{p}{q} \Bigr)^p |\nabla w_\alpha|^p \le (|u|+1)^{q-p} |\nabla u|^p \setone_{\{|u| \le \alpha\}}
			+ (\alpha+1)^{q-p} |\nabla u|^p \setone_{\{|u| > \alpha\}}.
	\end{equation}
	We will also need that
	\begin{equation}\label{eq:uglyalpha}
		\frac{1}{2} (\alpha+1)^{q-p}
			\le \frac{(\alpha+1)^{q-p+1} - 1}{\alpha + 1}
			\le (\alpha+1)^{q-p},
	\end{equation}
	which follows from the fact that $(\alpha+1)^{q-p+1} \ge 2$.

	From~\eqref{eq:nablavalpha}, \eqref{eq:nablawalpha}, \eqref{eq:uglyalpha} and~\eqref{eq:structure} we obtain that
	\begin{align*}
		& \int_\Omega \nabla v_\alpha \, A(x,u,\nabla u) \\
			& \quad \ge (q-p+1) \int_{\{|u| \le \alpha\}} (|u|+1)^{q-p} \Bigl( \nu |\nabla u|^p - \psi_1 |u|^p - \psi_1 \Bigr) \\
				& \qquad + \frac{(\alpha+1)^{q-p+1}-1}{\alpha+1} \int_{\{|u| > \alpha\}} \Bigl( \nu |\nabla u|^p - \psi_1 |u|^p - \psi_1 \Bigr) \\
			& \quad \ge \nu \int_{\{|u| \le \alpha\}} (|u|+1)^{q-p} |\nabla u|^p - 2 (q-p+1) \int_{\{|u| \le \alpha\}} \psi_1 (|u|+1)^q \\
				& \qquad + \frac{\nu}{2} \int_{\{|u| > \alpha\}} (\alpha+1)^{q-p} |\nabla u|^p - 2 (q-p+1) \int_{\{|u| > \alpha\}} \psi_1 (\alpha + 1)^{q-p} (|u|+1)^p \\
			& \quad \ge \frac{\nu}{2} \Bigl( \frac{p}{q} \Bigr)^p \int_\Omega |\nabla w_\alpha|^p - 2 (q-p+1) \int_\Omega \psi_1 w_\alpha^p.
	\end{align*}
	Similarly, we see that
	\begin{align*}
		& \Bigl| \int_\Omega v_\alpha \, B(x,u,\nabla u) \Bigr| \\
			& \le \int_{\{|u| \le \alpha\}} (|u|+1)^{q-p+1} \Bigl( \psi_3 |\nabla u|^{p-1} + \psi_1 |u|^{p-1} + \psi_1 \Bigr) \\
				& \qquad + \int_{\{|u| > \alpha\}} (\alpha+1)^{q-p} (|u| + 1) \Bigl( \psi_3 |\nabla u|^{p-1} + \psi_1 |u|^{p-1} + \psi_1 \Bigr) \\
			& \le \int_{\{|u| \le \alpha\}} \psi_3  \Bigl(\frac{q}{p}\Bigr)^{p-1} (|u|+1)^{\frac{q-p}{p}(p-1)} |\nabla u|^{p-1}  (|u|+1)^{\frac{q}{p}}
				+ 2 \int_{\{|u| \le \alpha\}} \psi_1 (|u|+1)^q \\
				& \qquad + \int_{\{|u| > \alpha\}} \psi_3 (\alpha + 1)^{\frac{q-p}{p}(p-1)} |\nabla u|^{p-1} (\alpha+1)^{\frac{q-p}{p}} (|u| + 1) \\
				& \qquad + 2 \int_{\{|u| > \alpha\}} \psi_1 (\alpha+1)^{q-p} (|u|+1)^p \\
			& = \int_\Omega \psi_3 |\nabla w_\alpha|^{p-1} w_\alpha + 2 \int_\Omega \psi_1 w_\alpha^p
	\end{align*}
	and
	\begin{equation*}
		\Bigl| \int_{\partial\Omega} v_\alpha h(x,u) \Bigr|
			\le 2 \int_{\partial\Omega} |v_\alpha| \; \psi_4 \; (|u|+1)^{p-1}
			\le 2 \int_{\partial\Omega} \psi_4 w_\alpha^p.
	\end{equation*}
	Using in addition that $u$ solves~\eqref{eq:RobinEqHom}, we have shown that
	\begin{align*}
		0 & = \int_\Omega \nabla v_\alpha \, A(x,u,\nabla u) + \int_\Omega v_\alpha B(x,u,\nabla u) + \int_{\partial\Omega} v_\alpha h(x,u) + \omega \int_\Omega v_\alpha u \\
			& \ge \frac{\nu}{2} \Bigl( \frac{p}{q} \Bigr)^p \int_\Omega |\nabla w_\alpha|^p - 2(q-p+2) \int_\Omega \psi_1 w_\alpha^p - \int_\Omega \psi_3 |\nabla w_\alpha|^{p-1} w_\alpha - 2 \int_{\partial\Omega} \psi_4 w_\alpha^p,
	\end{align*}
	i.e.,
	\begin{equation}\label{eq:walphagradest}
		\|w_\alpha\|_{W^{1,p}(\Omega)}^p \le c_q \int_\Omega \psi_1 w_\alpha^p + c_q \int_\Omega \psi_3 |\nabla w_\alpha|^{p-1} w_\alpha + c_q \int_{\partial\Omega} \psi_4 w_\alpha^p,
	\end{equation}
	where for simplicity we have assumed that $\psi_1 \ge 1$ almost everywhere, which constitutes no
	loss of generality.

	We now distinguish between the cases $p < N$, $p = N$ and $p > N$.
	\begin{enumerate}[(a)]
	\item
	Assume $p > N$. Then by~\cite[Theorem~5.8]{AF03} there exists $\theta \in (0,1)$ such that
	\[
		\|w_\alpha\|_{L^\infty(\partial\Omega)} \le \|w_\alpha\|_{L^\infty(\Omega)} \le c \|w_\alpha\|_{W^{1,p}(\Omega)}^\theta \|w_\alpha\|_{L^p(\Omega)}^{1-\theta}.
	\]
	Hence
	\begin{align*}
		& \int_\Omega \psi_1 w_\alpha^p + \int_\Omega \psi_3 |\nabla w_\alpha|^{p-1} w_\alpha + \int_{\partial\Omega} \psi_4 w_\alpha^p \\
			& \le \|\psi_1\|_{L^1(\Omega)} \|w_\alpha\|_{L^\infty(\Omega)}^p
				+ \|\psi_3\|_{L^p(\Omega)} \bigl\| |\nabla w_\alpha|^{p-1} \bigr\|_{L^{\frac{p}{p-1}}(\Omega)} \|w_\alpha\|_{L^\infty(\Omega)} \\
				& \qquad + \|\psi_4\|_{L^1(\partial\Omega)} \|w_\alpha\|_{L^\infty(\partial\Omega)}^p \\
			& \le c \|w_\alpha\|_{W^{1,p}(\Omega)}^{\theta p} \|w_\alpha\|_{L^p(\Omega)}^{(1-\theta)p}
				+ c \|w_\alpha\|_{W^{1,p}(\Omega)}^{p-1 + \theta} \|w_\alpha\|_{L^p(\Omega)}^{1-\theta} \\
			& \le \delta \|w_\alpha\|_{W^{1,p}(\Omega)}^p + c_\delta \|w_\alpha\|_{L^p(\Omega)}^p
	\end{align*}
	for every $\delta > 0$ by Young's inequality. Picking $\delta > 0$ small enough we obtain from this estimate
	and~\eqref{eq:walphagradest} that
	\[
		\|w_\alpha\|_{W^{1,p}(\Omega)}^p \le c_q \|w_\alpha\|_{L^p(\Omega)}^p.
	\]
	This proves~\eqref{eq:W1pest} since
	for $p = q$ we have $w_\alpha = |u|+1$ and $|\nabla w_\alpha| = |\nabla u|$, thus
	\[
		\|u\|_{W^{1,p}(\Omega)}^p
			\le \bigl\|w_\alpha\bigr\|_{W^{1,p}(\Omega)}^p
			\le c \|w_\alpha\|_{L^p(\Omega)}^p
			\le c \|u\|_{L^p(\Omega)}^p + c
	\]
	Finally, estimate~\eqref{eq:W1pest} implies~\eqref{eq:Lqest},
	in this particular case even for $q=\infty$,
	by the Sobolev embedding theorems.

	\item
	Now assume $p = N$. Then by the interpolation inequality for Lebesgue spaces and the Sobolev embedding theorems
	we find for every $r \in [1,\infty)$ an exponent $\theta_r \in (0,1)$ such that
	\[
		\|w_\alpha\|_{L^r(\Omega)}
			\le \|w_\alpha\|_{L^{2r}(\Omega)}^{\theta_r} \|w_\alpha\|_{L^1(\Omega)}^{1-\theta_r}
			\le c \|w_\alpha\|_{W^{1,p}(\Omega)}^{\theta_r} \|w_\alpha\|_{L^p(\Omega)}^{1-\theta_r}.
	\]
	Moreover, since the trace operator is compact from $W^{1,p}(\Omega)$ to $L^r(\partial\Omega)$ for every $r \in [1,\infty)$,
	we can estimate
	\[
		\|w_\alpha\|_{L^r(\partial\Omega)} \le \delta \|w_\alpha\|_{W^{1,p}(\Omega)} + c_\delta \|w_\alpha\|_{L^p(\Omega)}
	\]
	for every $\delta > 0$ by Lemma~\ref{lem:eberlein}.
	Using these two estimates and Young's inequality, we obtain with certain exponents $\theta$ and $\theta'$ in $(0,1)$ that
	\begin{align*}
		& \int_\Omega \psi_1 w_\alpha^p + \int_\Omega \psi_3 |\nabla w_\alpha|^{p-1} w_\alpha + \int_{\partial\Omega} \psi_4 w_\alpha^p \\
			& \le \|\psi_1\|_{L^{\frac{p}{p-\eps}}(\Omega)} \|w_\alpha^p\|_{L^{\frac{p}{\eps}}(\Omega)}
				+ \|\psi_3\|_{L^{\frac{p}{1-\eps}}(\Omega)} \bigl\| |\nabla w_\alpha|^{p-1} \bigr\|_{L^{\frac{p}{p-1}}(\Omega)} \|w_\alpha\|_{L^{\frac{p}{\eps}}(\Omega)} \\
				& \qquad + \|\psi_4\|_{L^{\frac{p-1}{p-1-\eps}}(\partial\Omega)} \|w_\alpha^p\|_{L^{\frac{p-1}{\eps}}(\partial\Omega)} \\
			& \le c \|w_\alpha\|_{W^{1,p}(\Omega)}^{p \theta} \|w_\alpha\|_{L^p(\Omega)}^{p (1-\theta)}
				+ c \|w_\alpha\|_{W^{1,p}(\Omega)}^{p-1 + \theta'} \|w_\alpha\|_{L^p(\Omega)}^{1-\theta'}
				+ c \|w_\alpha\|_{L^\frac{p(p-1)}{\eps}(\partial\Omega)}^p \\
			& \le \delta \|w_\alpha\|_{W^{1,p}(\Omega)}^p + c_\delta \|w_\alpha\|_{L^p(\Omega)}^p
	\end{align*}
	for every $\delta > 0$. As in the previous case $p > N$ this implies~\eqref{eq:W1pest} and
	hence~\eqref{eq:Lqest} by the Sobolev embedding theorems.

	\item
	Finally, assume $p < N$. Then
	\begin{align*}
		\int_\Omega \psi_1 w_\alpha^p
			& \le \|\psi_1\|_{L^{\frac{N}{p-\eps}}(\Omega)} \bigl\|w_\alpha^{p-\eps}\bigr\|_{L^{\frac{Np}{(N-p)(p-\eps)}}(\Omega)} \bigl\|w_\alpha^\eps\bigr\|_{L^{\frac{p}{\eps}}(\Omega)} \\
			& = c \|w_\alpha\|_{L^{\frac{Np}{N-p}}(\Omega)}^{p-\eps} \|w_\alpha\|_{L^p(\Omega)}^\eps
			\le c \|w_\alpha\|_{W^{1,p}(\Omega)}^{p-\eps} \|w_\alpha\|_{L^p(\Omega)}^\eps \\
			& \le \delta \|w_\alpha\|_{W^{1,p}(\Omega)}^p + c_\delta \|w_\alpha\|_{L^p(\Omega)}^p
	\end{align*}
	for every $\delta > 0$ by Young's inequality. Similarly,
	\begin{align*}
		& \int_\Omega \psi_3 |\nabla w_\alpha|^{p-1} w_\alpha \\
			& \qquad \le \|\psi\|_{L^{\frac{N}{1-\eps}}(\Omega)} \bigl\| |\nabla w_\alpha|^{p-1} \bigr\|_{L^{\frac{p}{p-1}}(\Omega)} \|w_\alpha^{1-\eps}\|_{L^{\frac{Np}{(N-p)(1-\eps)}}(\Omega)} \|w_\alpha^\eps\|_{L^{\frac{p}{\eps}}(\Omega)} \\
			& \qquad \le c \| w_\alpha \|_{W^{1,p}(\Omega)}^{p-1} \|w_\alpha\|_{L^{\frac{Np}{N-p}}(\Omega)}^{1-\eps} \|w_\alpha\|_{L^p(\Omega)}^\eps
			\le c \|w_\alpha\|_{W^{1,p}(\Omega)}^{p-\eps} \|w_\alpha\|_{L^p(\Omega)}^\eps \\
			& \qquad \le \delta \|w_\alpha\|_{W^{1,p}(\Omega)}^p + c_\delta \|w_\alpha\|_{L^p(\Omega)}^p.
	\end{align*}
	Moreover, since the trace operator is compact from $W^{1,p}(\Omega)$ to $L^{\frac{(N-1)p}{N-p+\eps}}(\partial\Omega)$,
	we can estimate
	\[
		\|w_\alpha\|_{L^{\frac{(N-1)p}{N-p+\eps}}(\partial\Omega)} \le \delta \|w_\alpha\|_{W^{1,p}(\Omega)} + c_\delta \|w_\alpha\|_{L^p(\Omega)}
	\]
	for every $\delta > 0$ by Lemma~\ref{lem:eberlein}.
	Thus, picking $\delta > 0$ small enough in the previous estimates,
	inequality~\eqref{eq:walphagradest} yields
	\begin{equation}\label{eq:walphaest}
		\|w_\alpha\|_{W^{1,p}(\Omega)}^p \le c_q \|w_\alpha\|_{L^p(\Omega)}^p.
	\end{equation}
	As in the previous cases this implies~\eqref{eq:W1pest}.

	Moreover, from~\eqref{eq:walphaest} and the Sobolev embedding theorems we deduce that
	\[
		\|w_\alpha\|_{L^{\frac{Np}{N-p}}(\Omega)}^p + \|w_\alpha\|_{L^{\frac{(N-1)p}{N-p}}(\partial\Omega)}^p
			\le c \|w_\alpha\|_{W^{1,p}(\Omega)}^p
			\le c_q \|w_\alpha\|_{L^p(\Omega)}^p.
	\]
	Since $w_\alpha$ approaches $(|u|+1)^{q/p}$ from below as $\alpha \to \infty$,
	the monotone convergence theorem implies that
	\[
		\bigl\| (|u|+1)^{q/p} \|_{L^{\frac{Np}{N-p}}(\Omega)}^p + \bigl\| (|u|+1)^{q/p} \bigr\|_{L^{\frac{(N-1)p}{N-p}}(\partial\Omega)}^p
			\le c_q \bigl\| (|u|+1)^{q/p} \bigr\|_{L^p(\Omega)},
	\]
	from which we can deduce that
	\begin{equation}\label{eq:uqest}
		\|u\|_{L^{\frac{Nq}{N-p}}(\Omega)} + \|u\|_{L^{\frac{(N-1)q}{N-p}}(\partial\Omega)}
			\le c_q \|u\|_{L^q(\Omega)} + c_q.
	\end{equation}
	Since this estimate holds for all $q \in [1,\infty)$, defining
	\[
		q_n \coloneqq \Bigl( \frac{N}{N-p} \Bigr)^n p
	\]
	and iterating~\eqref{eq:uqest} we obtain that
	\begin{align*}
		\|u\|_{L^{q_n}(\Omega)} + \|u\|_{L^{\frac{N-1}{N}q_n}(\partial\Omega)}
			& = \|u\|_{L^{\frac{Nq_{n-1}}{N-p}}(\Omega)} + \|u\|_{L^{\frac{(N-1)q_{n-1}}{N-p}}(\partial\Omega)} \\
			& \le c_{q_{n-1}} \|u\|_{L^{q_{n-1}}(\Omega)} + c_{q_{n-1}} \\
			& \le c_{q_{n-1}} c_{q_{n-2}} \|u\|_{L^{q_{n-2}}(\Omega)} + c_{q_{n-1}} c_{q_{n-2}} + c_{q_{n-1}} \\
			& \le \cdots \le c_n \|u\|_{L^{q_0}(\Omega)} + c_n
			= c_n \|u\|_{L^p(\Omega)} + c_n.
	\end{align*}
	Since $q_n \to \infty$ as $n \to \infty$, this proves~\eqref{eq:Lqest}.
	\qedhere
	\end{enumerate}
\end{proof}

\begin{remark}
	In the proof of Proposition~\ref{prop:Lqest} we silently passed over a subtlety that arises when deriving~\eqref{eq:walphagradest}.
	A priori we can test the equation~\eqref{eq:RobinEqHom} only against functions in $W^{1,p}(\Omega) \cap L^2(\Omega)$.
	However, we had to test the equation against $v_\alpha \in W^{1,p}(\Omega)$, but did not know a priori whether $v_\alpha \in L^2(\Omega)$.
	Still, since $\sgn(v_\alpha) = \sgn(u)$ we can pick a sequence $(\eta_n)$ of functions in $W^{1,p}(\Omega) \cap L^2(\Omega)$ that converges
	to $v_\alpha$ in $W^{1,p}(\Omega)$ and such that the sequence $(\eta_n \, u)$ is pointwise monotonically increasing.
	Then it follows from the monotone convergence theorem that
	\begin{align*}
		0 & = \int_\Omega \nabla\eta_n \, A(x,u,\nabla u) + \int_\Omega \eta_n \, B(x,u,\nabla u) + \omega \int_\Omega \eta_n \, u + \int_{\partial\Omega} \eta_n \, h(x,u) \\
			& \to \int_\Omega \nabla v_\alpha \, A(x,u,\nabla u) + \int_\Omega v_\alpha \, B(x,u,\nabla u) + \omega \int_\Omega v_\alpha \, u + \int_{\partial\Omega} v_\alpha \, h(x,u)
	\end{align*}
	as $n \to \infty$. Thus indeed equation~\eqref{eq:RobinEqHom} can be tested against $v_\alpha$, even before we know
	that $v_\alpha \in L^2(\Omega)$.
\end{remark}

Now that we have $L^q$-bounds at our disposal, it is easy to deduce the H\"older continuity of solutions from
the corresponding result concerning Neumann boundary conditions.
\begin{theorem}\label{thm:HoelderRobin}
	Let $\Omega$ be a Lipschitz domain and assume~\eqref{eq:structure} and~\eqref{eq:fFass}.
	Then there exists $\alpha \in (0,1)$ such that every weak solution $u \in W^{1,p}(\Omega)$ of~\eqref{eq:RobinEq} with $\omega \ge 0$
	is in $\mathrm{C}^{0,\alpha}(\Omega)$. Moreover, if $f$, $F$ and $g$ range over norm-bounded sets in their respective spaces, then
	the $\mathrm{C}^{0,\alpha}$-norms of the corresponding solutions $u$ remain bounded provided their norms in $L^p(\Omega)$
	remain bounded.
\end{theorem}
\begin{proof}
	Let $u \in W^{1,p}(\Omega)$ be a weak solution of~\eqref{eq:RobinEq}.
	Then $u$ is a weak solution of
	\begin{equation}\label{eq:HoelderRobin:mod}
		\left\{ \begin{aligned}
			-\divi \hat{A}(x,u,\nabla u) + \hat{B}(x,u,\nabla u) + \omega u & = 0 && \text{on } \Omega \\
			\hat{A}(x,u,\nabla u) \cdot \nu + \hat{h}(x,u) & = 0 && \text{on } \partial\Omega
		\end{aligned} \right.
	\end{equation}
	for $\hat{A}(x,u,z) \coloneqq A(x,u,z) - F(x)$, $\hat{B}(x,u,z) \coloneqq B(x,u,z) - f(x)$
	and $\hat{h}(x,u) \coloneqq h(x,u) - g(x)$. Assumption~\eqref{eq:fFass} ensures $\hat{A}$,
	$\hat{B}$ and $\hat{h}$ satisfy~\eqref{eq:structure}.

	If $p > N$, then $u \in \mathrm{C}^{0,\alpha}(\Omega)$ by a Sobolev embedding theorem.
	The boundedness assertion about $u$ follows from Proposition~\ref{prop:Lqest} applied to the
	equation~\eqref{eq:HoelderRobin:mod} since the 
	constants in~\eqref{eq:W1pest} depend only on upper bounds for the structure parameters of $\hat{A}$,
	$\hat{B}$ and $\hat{h}$, which in turn depend only on $A$, $B$, $h$ and upper bound for $f$, $F$ and $g$.

	Now assume $p \le N$.
	Pick $q \in [p,\infty)$ such that $q \ge \frac{2(N-1)(p-1)}{\eps}$ and $q \ge \frac{N}{p-\eps}$
	with $\eps > 0$ as in the structure condition~\eqref{eq:structure}.
	We obtain from Proposition~\ref{prop:Lqest} applied to equation~\eqref{eq:HoelderRobin:mod}
	that $u \in L^q(\Omega)$ and $u|_{\partial\Omega} \in L^q(\partial\Omega)$ with a bound that depends
	only on bounds for the structure parameters of $\hat{A}$, $\hat{B}$ and $\hat{h}$. Thus
	\[
		h_0 \coloneqq h(\cdot,u(\cdot)) \in L^{\frac{N-1}{p-1-\eps/2}}(\partial\Omega)
	\]
	by the structure assumption~\eqref{eq:structure}.
	Since $u$ is a weak solution of the Neumann problem
	\begin{equation}\label{eq:HoelderRobin:Neumod}
		\left\{ \begin{aligned}
			-\divi A(x,u,\nabla u) + B(x,u,\nabla u) & = f(x) - \omega u(x) - (\divi F)(x) && \text{on } \Omega \\
			A(x,u,\nabla u) \cdot \nu & = g(x) - h_0(x) + (F \cdot \nu)(x) && \text{on } \partial\Omega,
		\end{aligned} \right.
	\end{equation}
	we obtain that $u \in \mathrm{C}^{0,\alpha}(\Omega)$ with the same $\alpha \in (0,1)$ as in
	Theorem~\ref{thm:HoelderNeumann}. The boundedness assertion about $u$ in $\mathrm{C}^{0,\alpha}(\Omega)$
	follows from~\eqref{eq:HoelderNeumann:est} and the observation that the right hand side of~\eqref{eq:HoelderRobin:Neumod}
	can be estimated in terms of $A$, $B$, $h$, $\omega$ and upper bounds for $f$, $F$ and $h$.
\end{proof}

\begin{remark}
	Working with the full equation~\eqref{eq:RobinEq} instead of~\eqref{eq:RobinEqHom}
	in Proposition~\ref{prop:Lqest}, we could have found an estimate like~\eqref{eq:HoelderNeumann:est}
	also in the situation of Theorem~\eqref{thm:HoelderRobin}.
	However, for our purposes nothing is gained by this more precise estimate, so for the sake of simplicity
	we avoided this additional effort.
\end{remark}

Sometimes it is convenient to replace the functions on the right hand side of~\eqref{eq:RobinEq}
by an arbitrary bounded linear functional on $W^{1,p}(\Omega)$. Therefore we formulate the regularity result
of Theorem~\ref{thm:HoelderRobin} also for such equations.

\begin{corollary}
	Let $\Omega \subset \mathds{R}^N$ be a Lipschitz domain and assume~\eqref{eq:structure}.
	Let $q \in [1,\infty)$ satisfy
	\[
		\left\{ \begin{aligned}
			q & \le \frac{N}{N-p+1} && \text{if } p < N \\
			q & < N && \text{if } p = N \\
			q & \le p && \text{if } p > N
		\end{aligned} \right.
	\]
	and fix $\psi \in W^{1,q}(\Omega)'$. Then every function $u \in W^{1,p}(\Omega)$ that satisfies
	\[
		\int_\Omega \nabla\eta \, A(x,u,\nabla u) + \int_\Omega \eta \, B(x,u,\nabla u) + \int_{\partial\Omega} \eta \, h(x,u) = \psi(\eta)
	\]
	for all $\eta \in \mathrm{C}^\infty(\overline{\Omega})$ is H\"older continuous on $\Omega$.
\end{corollary}
\begin{proof}
	This is a direct consequence of the representation formula for functionals on Sobolev spaces~\cite[\S 4.3]{Ziemer}
	and Theorem~\ref{thm:HoelderRobin}.
\end{proof}

\section{Parabolic problems}\label{sec:par}

Let $\Omega \subset \mathds{R}^N$ be a Lipschitz domain. We show that the parabolic problem
\begin{equation}\label{eq:ParRobin}
	\left\{ \begin{aligned}
		u_t(t,x) - \divi a(x,\nabla u(t,x)) + b(x,u(t,x)) & = 0 && t > 0, \; x \in \Omega \\
		a(x,\nabla u(t,x)) \cdot \nu + h(x,u(x)) & = 0 && t > 0, \; x \in \partial\Omega \\
		u(0,x) & = u_0(x) && x \in \Omega
	\end{aligned} \right.
\end{equation}
with Robin boundary conditions is well-posed in the space $\mathrm{C}(\overline{\Omega})$ under
suitable conditions on $a$, $b$ and $h$.
More precisely, we assume that $a\colon \Omega \times \mathds{R}^N \to \mathds{R}^N$,
$b\colon \Omega \times \mathds{R} \to \mathds{R}$
and $h\colon \partial\Omega \times \mathds{R} \to \mathds{R}$ are measurable
and continuous in the second argument.
We also require that there exists $p \in (1,\infty)$ such that
\begin{equation}\label{eq:structure2}
	\left\{ \begin{aligned}
		z \, a(x,z) & \ge \nu |z|^p - \psi_1(x), &
		|a(x,z)| & \le \mu |z|^{p-1} + \psi_2(x), \\
		|b(x,u)| & \le \psi_1(x) |u|^{p-1} + \psi_1(x), &
		|h(x,u)| & \le \psi_4(x) |u|^{p-1} + \psi_4(x)
	\end{aligned} \right.
\end{equation}
for all $x \in \Omega$, $u \in \mathds{R}$ and $z \in \mathds{R}^N$
with $\psi_1$, $\psi_2$ and $\psi_4$ as in~\eqref{eq:structure},
i.e., the functions $A(x,u,z) \coloneqq a(x,z)$, $B(x,u,z) \coloneqq b(x,z)$ and $h$ satisfy~\eqref{eq:structure}.
Moreover, we assume the monotonicity conditions
\begin{equation}\label{eq:monass}
	\left\{ \begin{aligned}
		(z_1 - z_2) \, (a(x,z_1) - a(x,z_2)) & \ge 0, \\
		(u_1 - u_2) \, (b(x,u_1) - b(x,u_2)) & \ge 0, \\
		(u_1 - u_2) \, (h(x,u_1) - h(x,u_2)) & \ge 0
	\end{aligned} \right.
\end{equation}
for all $x \in \Omega$, $u_1,u_2 \in \mathds{R}$ and $z_1,z_2 \in \mathds{R}^N$.

In order to prove well-posedness in $\mathrm{C}(\overline{\Omega})$ we are
going to show that the operator which is naturally associated
with~\eqref{eq:ParRobin} is m-accretive on $\mathrm{C}(\overline{\Omega})$ and
thus generates a semigroup on that space.

It is convenient to first introduce a version $\mathcal{A}$ of the operator
associated with~\eqref{eq:ParRobin} acting from $V \coloneqq W^{1,p}(\Omega)
\cap L^2(\Omega)$ into its $V'$ and to study its properties. Later on we will
turn our attention to its part in $\mathrm{C}(\overline{\Omega})$ in order to
obtain the main result.
\begin{proposition}\label{prop:monop}
	The space $V \coloneqq W^{1,p}(\Omega) \cap L^2(\Omega)$ is a separable, reflexive Banach space for the norm
	$\|u\|_V \coloneqq \|u\|_{W^{1,p}(\Omega)} + \|u\|_{L^2(\Omega)}$. Under the assumptions~\eqref{eq:structure2}
	and~\eqref{eq:monass}, via
	\begin{equation}\label{eq:cAdef}
		\apply{\mathcal{A}u}{v} \coloneqq \int_\Omega \nabla v \; a(x,\nabla u) + \int_\Omega v \, b(x,u) + \int_{\partial\Omega} v \, h(x,u).
	\end{equation}
	we have defined a bounded, continuous, monotone operator $\mathcal{A}\colon V \to V'$.
\end{proposition}
\begin{proof}
	We prove the assertions only for $p < N$ and only mention that the case $p \ge N$ is similar.
	Identifying $V$ with a closed subspace of the direct sum $W^{1,p}(\Omega) \oplus L^2(\Omega)$, we have proved the first claim.
	For $u$ and $v$ in $W^{1,p}(\Omega)$ we have by~\eqref{eq:structure2} that
	\begin{align*}
		\int_\Omega \bigl| \nabla v \; a(x,\nabla u) \bigr|
			& \le \|\nabla v\|_{L^p(\Omega)} \Bigl( \int_\Omega \bigl( \mu |\nabla u|^{p-1} + \psi_2 \bigr)^{\frac{p}{p-1}} \Bigr)^{\frac{p-1}{p}} \\
			& \le c \|v\|_{W^{1,p}(\Omega)} \bigl( \|\nabla u\|_{L^p(\Omega)}^{p-1} + \|\psi_2\|_{L^{\frac{p}{p-1}}(\Omega)} \bigr) \\
			& \le \bigl( c \|u\|_V^{p-1} + c \bigr) \|v\|_V.
	\end{align*}
	Similarly, using in addition the Sobolev embedding theorems, we obtain that
	\begin{align*}
		 \int_\Omega \bigl| v \, b(x,u) \bigr|
			& \le \|v\|_{L^{\frac{Np}{N-p}}(\Omega)} \|\psi_1\|_{L^{\frac{N}{p}}(\Omega)} \bigl\| |u|^{p-1} + 1 \bigr\|_{L^{\frac{Np}{(N-p)(p-1)}}(\Omega)} \\
			& \le c \|v\|_{W^{1,p}(\Omega)} \bigl( c \|u\|_{W^{1,p}(\Omega)}^{p-1} + c \bigr)
			\le \bigl( c \|u\|_V^{p-1} + c \bigr) \|v\|_V
	\end{align*}
	and
	\begin{align*}
		 \int_{\partial\Omega} \bigl| v \, h(x,u) \bigr|
			& \le \|v\|_{L^{\frac{(N-1)p}{N-p}}(\partial\Omega)} \|\psi_4\|_{L^{\frac{N-1}{p-1}}(\partial\Omega)} \bigl\| |u|^{p-1} + 1 \bigr\|_{L^{\frac{(N-1)p}{(N-p)(p-1)}}(\partial\Omega)} \\
			& \le c \|v\|_{W^{1,p}(\Omega)} \bigl( c \|u\|_{W^{1,p}(\Omega)}^{p-1} + c \bigr)
			\le \bigl( c \|u\|_V^{p-1} + c \bigr) \|v\|_V.
	\end{align*}
	Thus $\mathcal{A}$ is well-defined and $\|\mathcal{A}u\|_{V'} \le c \|u\|_V^{p-1} + c$ for all $u \in V$,
	proving that the operator $\mathcal{A}$ is bounded, i.e., that $\mathcal{A}$ maps bounded sets into bounded sets.

	Now we show that $\mathcal{A}$ is continuous.
	To this end, let $(u_n)$ be a sequence in $V$ that converges to $u \in V$. Passing to a subsequence,
	we can assume that $\nabla u_n \to \nabla u$ pointwise and that $|\nabla u_n| \le m$ for some function $m \in L^p(\Omega)$.
	Then $a(x,\nabla u_n(x)) \to a(x,\nabla u(x))$ for almost every $x \in \Omega$ and
	\[
		|a(x,\nabla u_n)| \le \mu m^{p-1} + \psi_2 \in L^{\frac{p}{p-1}}(\Omega).
	\]
	Hence $a(\cdot,\nabla u_n) \to a(\cdot,\nabla u)$ in $L^{\frac{p}{p-1}}(\Omega)$ by the dominated convergence theorem.
	An analogous reasoning yields that
	\[
		b(\cdot,u_n) \to b(\cdot,u) \text{ in } L^{\frac{Np}{Np-N+p}}(\Omega)
		\quad\text{and}\quad
		h(\cdot,u_n) \to h(\cdot,u) \text{ in } L^{\frac{(N-1)p}{N(p-1)}}(\partial\Omega).
	\]
	Thus, given $\eps > 0$, we have
	\begin{align*}
		|\apply{\mathcal{A}u_n - \mathcal{A}u}{v}|
			& \le \|\nabla v\|_{L^p(\Omega)} \|a(\cdot,\nabla u_n) - a(\cdot,\nabla u)\|_{L^{\frac{p}{p-1}}(\Omega)} \\
				& \qquad + \|v\|_{L^{\frac{Np}{N-p}}(\Omega)} \|b(\cdot,u_n) - b(\cdot,u)\|_{L^{\frac{Np}{Np-N+p}}(\Omega)} \\
				& \qquad + \|v\|_{L^{\frac{(N-1)p}{N-p}}(\partial\Omega)} \|h(\cdot,u_n) - h(\cdot,u)\|_{L^{\frac{(N-1)p}{N(p-1)}}(\partial\Omega)} \\
			& \le \eps \|v\|_{W^{1,p}(\Omega)}
	\end{align*}
	for $n \ge n_0(\eps)$, which shows that $\mathcal{A}u_n \to \mathcal{A}u$ in $V'$.

	Finally, the monotonicity of $\mathcal{A}$, i.e., $\apply{\mathcal{A}u - \mathcal{A}v}{u-v} \ge 0$ for all $u,v \in V$,
	is a trivial consequence of~\eqref{eq:monass}.
\end{proof}

Next we show that $\mathcal{A}$ is bijective.
\begin{proposition}\label{prop:maxmon}
	Under the assumptions of Proposition~\ref{prop:monop}, for every $\phi \in V'$ and every $\alpha > 0$
	there exists a unique function $u \in V$ such that
	\begin{equation}\label{eq:maxmon:eq}
		\scalar[L^2(\Omega)]{u}{v} + \alpha \apply{\mathcal{A}u}{v} = \phi(v)
	\end{equation}
	for all $v \in V$.
\end{proposition}
\begin{proof}
	Define the operator $\mathcal{A}_\alpha\colon V \to V'$ by
	\[
		\apply{\mathcal{A}_\alpha u}{v} \coloneqq \scalar[L^2(\Omega)]{u}{v} + \alpha \apply{\mathcal{A}u}{v}.
	\]
	From Proposition~\ref{prop:monop} we obtain that $\mathcal{A}_\alpha$ is bounded, continuous and monotone.
	From~\eqref{eq:structure2} we obtain that
	\[
		\int_\Omega \nabla u \; a(x,\nabla u)
			\ge \nu \int_\Omega |\nabla u|^p - \int_\Omega \psi_1
			= \nu \|\nabla u\|_{L^p(\Omega)}^p - c.
	\]
	for all $u \in V$.
	Moreover, \eqref{eq:structure2} and~\eqref{eq:monass} yield
	\begin{align*}
		\int_\Omega u \; b(x,u)
			& = \int_\Omega (u-0) \, \bigl( b(x,u) - b(x,0) \bigr) + \int_\Omega u \, b(x,0)
			\ge -\int_\Omega \psi_1 |u| \\
			& \ge - \|\psi_1\|_{L^{\frac{N}{p}}(\Omega)} \|u\|_{L^{\frac{N}{N-p}}(\Omega)}
			\ge -c \|u\|_{W^{1,p}(\Omega)}
	\end{align*}
	for all $u \in V$, and analogously we see that $\int_\Omega u \, h(x,u) \ge -c \|u\|_{W^{1,p}(\Omega)}$.
	Combining the latter three estimates we have shown that
	\begin{equation}\label{eq:coerc}
		\begin{aligned}
			\apply{\mathcal{A}_\alpha u}{u}
				& \ge \|u\|_{L^2(\Omega)}^2 + \alpha \nu \|\nabla u\|_{L^p(\Omega)}^p - c \|u\|_{W^{1,p}(\Omega)} \\
				& \ge \beta \bigl( \|u\|_{L^2(\Omega)} + \|\nabla u\|_{L^p(\Omega)} \bigr)^q - c \|u\|_V - c
		\end{aligned}
	\end{equation}
	for $q \coloneqq \min\{2,p\} > 1$ and some $\beta > 0$.

	From Lemma~\eqref{lem:eberlein} we obtain that
	\[
		\|u\|_{L^p(\Omega)} \le \delta \|u\|_V + c_\delta \|u\|_{L^2(\Omega)}
			\le c \delta \|\nabla u\|_{L^p(\Omega)} + c \delta \|u\|_{L^p(\Omega)} + c_\delta \|u\|_{L^2(\Omega)}
	\]
	for all $u \in V$ with an arbitrary $\delta > 0$. Picking $\delta > 0$ small enough we deduce that
	\[
		\|u\|_{L^p(\Omega)} \le c \|\nabla u\|_{L^p(\Omega)} + c \|u\|_{L^2(\Omega)}.
	\]
	Hence $\|u\|_{2,p} \coloneqq \|u\|_{L^2(\Omega)} + \|\nabla u\|_{L^p(\Omega)}$ is an equivalent norm on $V$.
	Thus from~\eqref{eq:coerc} we obtain that
	\begin{equation}\label{eq:coercive}
		\liminf_{\|u\|_V \to \infty} \frac{\apply{\mathcal{A}_\alpha u}{u}}{\|u\|_V}
			\ge \liminf_{\|u\|_{2,p} \to \infty} \frac{ \beta \|u\|_{2,p}^q - c \|u\|_{2,p} - c }{ \|u\|_{2,p} }
			= \infty
	\end{equation}
	and call~\eqref{eq:coerc} the coercivity of $\mathcal{A}_\alpha$. Since $\mathcal{A}_\alpha$ is also bounded, continuous
	and monotone, the operator is surjective by the Minty-Browder theorem~\cite[\S II.2]{Showalter}.
	This means that~\eqref{eq:maxmon:eq} has a solution for every $\phi \in V'$.
	Moreover, if $\mathcal{A}_\alpha u_1 = \mathcal{A}_\alpha u_2$
	for two functions $u_1$ and $u_2$ in $V$, then by monotonicity of $\mathcal{A}$ we obtain that
	\[
		0 = \apply{\mathcal{A}_\alpha u_1 - \mathcal{A}_\alpha u_2}{u_1 - u_2}
			\ge \scalar[L^2(\Omega)]{u_1 - u_2}{u_1 - u_2} = \|u_1 - u_2\|_{L^2(\Omega)}^2,
	\]
	implying that $u_1 = u_2$. Hence the solution of~\eqref{eq:maxmon:eq} is unique.
\end{proof}

Our next step is to show that $(I + \alpha \mathcal{A})^{-1}$ is a contraction with respect to the norms
of $L^q(\Omega)$ for $q \in [2,\infty]$. This is equivalent to saying that the part $\mathcal{A}$ in $L^q(\Omega)$ is
accretive. Setting $q = \infty$ we thus obtain in particular that the part of $\mathcal{A}$ in $\mathrm{C}(\overline{\Omega})$ is accretive.
\begin{proposition}\label{prop:maxacc}
	Let~\eqref{eq:structure2} and~\eqref{eq:monass} be satisfied. Let $f_1$ and $f_2$ be in $L^q(\Omega) \subset L^2(\Omega)$,
	$q \in [2,\infty]$ and define $\phi_i(v) \coloneqq \scalar[L^2(\Omega)]{f_i}{v}$ so that $\phi_i \in V'$.
	Then the unique solutions $u_i$ of~\eqref{eq:maxmon:eq}
	for the right hand sides $\phi_i$ satisfy $\|u_1 - u_2\|_{L^q(\Omega)} \le \|f_1 - f_2\|_{L^q(\Omega)}$.
	Moreover, if $q > \frac{N}{p}$, then $u_1$ and $u_2$ are in $\mathrm{C}(\overline{\Omega})$.
\end{proposition}
\begin{proof}
	Fix $q \in [2,\infty)$ and $k \ge 1$. Then
	by the chain rule~\cite[Theorem~2.2.2]{Ziemer} the function
	\[
		v_k \coloneqq \bigl( |u_1 - u_2| \wedge k \bigr)^{q-2} (u_1 - u_2),
	\]
	where $x \wedge y$ denotes the minimum of $x$ and $y$, is in $V$ with weak derivative
	\[
		\nabla v_k = (q-1) |u_1 - u_2|^{q-2} (\nabla u_1 - \nabla u_2) \setone_{\{|u_1 - u_2| \le k\}} + k^{q-2} (\nabla u_1 - \nabla u_2) \setone_{\{|u_1 - u_2| > k\}}.
	\]
	Since at each point $v_k$ is a positive multiple of $u_1-u_2$ and $\nabla v_k$ is a positive multiple of $\nabla u_1 - \nabla u_2$,
	we deduce from~\eqref{eq:monass} that
	\[
		\apply{\mathcal{A}u_1 - \mathcal{A}u_2}{v_k} \ge 0.
	\]
	Hence with $\mathcal{A}_\alpha$ as in the proof of Proposition~\ref{prop:maxmon} we obtain that
	\begin{align*}
		& \int_\Omega \bigl( |u_1 - u_2| \wedge k \bigr)^{q-2} |u_1 - u_2|^2
			= \scalar[L^2(\Omega)]{u_1-u_2}{v_k}
			\le \apply{\mathcal{A}_\alpha u_1 - \mathcal{A}_\alpha u_2}{v_k} \\
			& \qquad = \phi_1(v_k) - \phi_2(v_k)
			\le \|f_1 - f_2\|_{L^q(\Omega)} \|v_k\|_{L^{\frac{q}{q-1}}(\Omega)} \\
			& \qquad \le \|f_1 - f_2\|_{L^q(\Omega)} \Bigl( \int_\Omega \bigl( |u_1 - u_2| \wedge k \bigr)^{q-2} |u_1 - u_2|^2 \Bigr)^{\frac{q-1}{q}}
	\end{align*}
	Dividing in this equation and afterwards letting $k$ tend to infinity, we obtain from the monotone convergence theorem that
	$\|u_1 - u_2\|_{L^q(\Omega)} \le \|f_1 - f_2\|_{L^q(\Omega)}$ for every $q \in [2,\infty)$. If $f_1$ and $f_2$ are in $L^\infty(\Omega)$,
	we pass to the limit $q \to \infty$ and obtain that $\|u_1 - u_2\|_{L^\infty(\Omega)} \le \|f_1 - f_2\|_{L^\infty(\Omega)}$.

	Now assume $q > \frac{N}{p}$.
	By definition, the function $u_i$ is a weak solution of
	\[
		\left\{ \begin{aligned}
			-\divi a(x, \nabla u) + b(x, u) & = f_i && \text{on } \Omega \\
			a(x,\nabla u) \cdot \nu + h(x,u) & = 0 && \text{on } \partial\Omega.
		\end{aligned} \right.
	\]
	Hence $u_i \in \mathrm{C}^{0,\alpha}(\Omega) \subset \mathrm{C}(\overline{\Omega})$ by Theorem~\ref{thm:HoelderRobin}.
\end{proof}

Regard $\mathrm{C}(\overline{\Omega})$ as a subspace of $V'$ by identifying
a function $f \in \mathrm{C}(\overline{\Omega})$ with the linear functional $v \mapsto \scalar[L^2(\Omega)]{f}{v}$.
Then the part $\mathcal{A}^c$ of $\mathcal{A}$ in $\mathrm{C}(\overline{\Omega})$ is the restriction of $\mathcal{A}$ to
\[
	D(\mathcal{A}^c) \coloneqq \bigl\{ u \in V \cap \mathrm{C}(\overline{\Omega}) : \mathcal{A}u \in \mathrm{C}(\overline{\Omega}) \bigr\}.
\]
We regard $\mathcal{A}^c$ as a non-linear (and single-valued) operator on $\mathrm{C}(\overline{\Omega})$.
Proposition~\ref{prop:maxacc} implies that $\mathcal{A}^c$ is m-accretive,
i.e., for all $\alpha > 0$ the operator $I + \alpha \mathcal{A}^c\colon D(\mathcal{A}^c) \to \mathrm{C}(\overline{\Omega})$
is bijective and $(I+\alpha\mathcal{A}^c)^{-1}$ is a contraction on $\mathrm{C}(\overline{\Omega})$.
We now show that $\mathcal{A}^c$ is densely defined.

\begin{proposition}\label{prop:dense}
	Under assumptions~\eqref{eq:structure2} and~\eqref{eq:monass} the set $D(\mathcal{A}^c)$ is dense in $\mathrm{C}(\overline{\Omega})$.
\end{proposition}
\begin{proof}
	We give the proof only for $p < N$ and only mention that the case $p \ge N$ can be treated analogously.
	First assume that $u \in \mathrm{C}^\infty(\overline{\Omega})$, so in particular $u \in V$.
	Then $\phi \coloneqq \mathcal{A}_1 u \in V'$, where
	$\mathcal{A}_1$ is defined as in the proof of Proposition~\ref{prop:maxmon}. More precisely,
	\begin{align*}
		|a(x,\nabla u)| & \le \mu \|\nabla u\|_{L^\infty(\Omega)}^{p-1} + \psi_2 \le c_u + \psi_2 \\
		|b(x,u)| & \le \psi_1 (\|u\|_{L^\infty(\Omega)}^{p-1} + 1) \le c_u \psi_1 \\
		|h(x,u)| & \le \psi_4 (\|u\|_{L^\infty(\Omega)}^{p-1} + 1) \le c_u \psi_4
	\end{align*}
	and hence
	\begin{align*}
		\phi(v)
			& = \int_\Omega vu + \int_\Omega \nabla v \, a(x,\nabla u) + \int_\Omega v \, b(x,u) + \int_{\partial\Omega} v \, h(x,u) \\
			& \le c_u \|v\|_{L^1(\Omega)}
				+ \|\nabla v\|_{L^{\frac{N}{N-p+1}}(\Omega)} \| c_u
				+ \psi_2 \|_{L^{\frac{N}{p-1}}(\Omega)} \\
				& \qquad + c_u \|v\|_{L^{\frac{N}{N-p}}(\Omega)} \|\psi_1\|_{L^{\frac{N}{p}}(\Omega)}
				+ c_u \|v\|_{L^{\frac{N-1}{N-p}}(\partial\Omega)} \|\psi_4\|_{L^{\frac{N-1}{p-1}}(\partial\Omega)} \\
			& \le c_u \|v\|_{W^{1,\frac{N}{N-p+1}}(\Omega)}
	\end{align*}
	for all $v \in V$. Thus $\psi$ extends to a bounded functional on $W^{1,\frac{N}{N-p+1}}(\Omega)$.
	Consequently, there exist $k \in L^{\frac{N}{p-1}}(\Omega)$ and $K \in L^{\frac{N}{p-1}}(\Omega;\mathds{R}^N)$
	such that
	\[
		\phi(v) = \int_\Omega v \, k + \int_\Omega \nabla v \; K
	\]
	for all $v \in V$, see~\cite[\S 4.3]{Ziemer}. Pick sequences $k_n \in \mathrm{C}^\infty_c(\Omega)$ and
	$K_n \in \mathrm{C}^\infty_c(\Omega;\mathds{R}^N)$ such that $k_n \to k$ and $K_n \to K$ in the $L^{\frac{N}{p-1}}$-norm.
	Then $f_n \coloneqq k_n - \divi K_n$ is in $\mathrm{C}^\infty_c(\Omega) \subset \mathrm{C}(\overline{\Omega})$.
	Thus by Propositions~\ref{prop:maxmon} and the additional claim in~\eqref{prop:maxacc}
	there exists $u_n \in D(\mathcal{A}^c)$ such that $(I + \mathcal{A}^c)u_n = f_n$.

	Define $\phi_n(v) \coloneqq \scalar[L^2(\Omega)]{f_n}{v}$, which can equivalently be written as $\phi_n \coloneqq \mathcal{A}_1 u_n$. Then
	\begin{align*}
		|\phi_n(v) - \phi(v)|
			& \le \int_\Omega |v| \, |k_n - k| + \int_\Omega |\nabla v| \, |K_n - K| \\
			& \le \|v\|_{W^{1,p}(\Omega)} \bigl( \|k_n - k\|_{L^{\frac{p}{p-1}}(\Omega)} + \|K_n - K\|_{L^{\frac{p}{p-1}}(\Omega)} \bigr).
	\end{align*}
	This shows that $\phi_n \to \phi$ in $V'$.
	In particular we see that $(\phi_n)$ is bounded in $V'$, which implies that
	\[
		\apply{\mathcal{A}_1 u_n}{u_n} = \phi_n(u_n) \le c_u \|u_n\|_V.
	\]
	By~\eqref{eq:coercive} this implies that $(u_n)$ is bounded in $V$. Thus passing to a subsequence
	we can assume that $(u_n)$ converges weakly to some $\tilde{u} \in V$.

	By Minty's theorem~\cite[Proposition~II.2.2]{Showalter} we have
	\[
		\apply{\mathcal{A}_1 v - \phi_n}{v - u_n} \ge 0
	\]
	for all $v \in V$.
	Since $\phi_n \to \phi$ in $V'$ and $u_n \wto \tilde{u}$ in $V$ we obtain by taking limits that
	\[
		\apply{\mathcal{A}_1 v - \phi}{v - \tilde{u}} \ge 0
	\]
	for all $v \in V$. Using Minty's theorem once again we deduce from this that
	$\mathcal{A}_1 \tilde{u} = \phi = \mathcal{A}_1 u$. By the uniqueness assertion of
	Proposition~\ref{prop:maxmon} this implies $\tilde{u} = u$. Thus we have shown that $u_n \wto u$ in $V$.

	We have seen that $(u_n)$ is bounded in $V$ and hence in particular in $L^p(\Omega)$.
	Since by construction $u_n$ is a weak solution of
	\[
		\left\{ \begin{aligned}
			-\divi a(x,\nabla u_n) + b(x,u_n) & = k_n - \divi K_n && \text{on } \Omega \\
			a(x,\nabla u_n) \cdot \nu + h(x,u_n) & = K_n \cdot \nu && \text{on } \partial\Omega,
		\end{aligned} \right.
	\]
	we obtain from Theorem~\ref{thm:HoelderRobin} that $(u_n)$
	is bounded in $\mathrm{C}^{0,\alpha}(\Omega)$. Exploiting compactness we can assume after passing to a subsequence that
	$(u_n)$ has a limit in $\mathrm{C}(\overline{\Omega})$. Since $u_n \wto u$ in $V$, this limit is $u$.
	Hence we have shown that for each $u \in \mathrm{C}^\infty(\overline{\Omega})$ there exists a sequence $u_n \in D(\mathcal{A}^c)$
	such that $u_n \to u$ in $\mathrm{C}(\overline{\Omega})$. This proves that $\mathrm{C}^\infty(\overline{\Omega})$
	is contained in the closure of $D(\mathcal{A}^c)$ in $\mathrm{C}(\overline{\Omega})$. Since $\mathrm{C}^\infty(\overline{\Omega})$
	is dense in $\mathrm{C}(\overline{\Omega})$, we have shown that $D(\mathcal{A}^c)$ is dense in $\mathrm{C}(\overline{\Omega})$.
\end{proof}

In order to state the main result of this section it is convenient to introduce the notion of a non-linear contraction $\mathrm{C}_0$-semigroup.
The definition is consistent with the linear case, i.e., a family of linear operators on a Banach space $X$ is a
non-linear contraction $\mathrm{C}_0$-semigroup if and only if it is a linear contraction $\mathrm{C}_0$-semigroup in the
usual sense of for example~\cite{HP57}.
\begin{definition}\label{def:semigroup}
	Let $X$ be a Banach space and
	let $\mathcal{B}\colon X \supset D(\mathcal{B}) \to X$ be an m-accretive operator on $X$,
	i.e., for all $\alpha > 0$ the operator $I + \alpha \mathcal{B}\colon D(\mathcal{B}) \to X$ is bijective
	with $(I + \alpha \mathcal{B})^{-1}\colon X \to X$ being a contraction.
	Then by the Crandall-Liggett theorem~\cite[\S IV.8]{Showalter} the limit
	$S(t)u_0 = \lim_{n \to \infty} (I + \frac{t}{n} \mathcal{B})^{-n} u_0$ exists for $u_0 \in C \coloneqq \overline{D(\mathcal{B})}$
	and $t \ge 0$ and the mappings $S(t)$ satisfy
	\begin{enumerate}[(i)]
	\item
		$S(t)\colon C \to C$ is contractive for every $t \ge 0$;
	\item
		$S(0) = \id_C$;
	\item
		$S(t+s) = S(t) \circ S(s)$ for all $t, s \ge 0$;
	\item
		$t \mapsto S(t)u_0$ is continuous for all $u_0 \in C$.
	\end{enumerate}
	We say that $(S(t))_{t \ge 0}$ is a non-linear contraction $\mathrm{C}_0$-semigroup on $C$
	and call $-\mathcal{B}$ its generator.
\end{definition}

The following remark about non-linear contraction $\mathrm{C}_0$-semigroups and their generators justifies
that we regard the trajectories of a non-linear contraction $\mathrm{C}_0$-semigroup with generator $-\mathcal{B}$ as
the unique solutions of the problem $u'(t) + \mathcal{B}u(t) = 0$.
For the definition of a $\mathrm{C}^0$-solution, which is frequently also called mild solution in the literature, and proofs of the following
fact we refer to~\cite[\S IV.3 and \S IV.8]{Showalter} or~\cite{BCP}.

\begin{remark}\label{rem:generator}
	If $\mathcal{B}$ is m-accretive and $(S(t))_{t \ge 0}$ is the semigroup generated by $-\mathcal{B}$,
	then for each $u_0 \in C \coloneqq \overline{D(\mathcal{B})}$ and $T > 0$ the unique $\mathrm{C}^0$-solution
	of $u'(t) + \mathcal{B}u(t) = 0$, $u(0) = u_0$ on $[0,T]$ is given by $u(t) = S(t)u_0$.

	Moreover, if $X$ is a Hilbert space and $u_0 \in D(\mathcal{B})$, then this unique solution $u$ is Lipschitz continuous,
	$u(t) \in D(\mathcal{B})$ for all $t \ge 0$ and $u'(t) + \mathcal{B}u(t) = 0$ for almost every $t > 0$.
	In this situation we say that $u$ is a strong solution of $u'(t) + \mathcal{B}u(t) = 0$.
\end{remark}

Now regard $L^2(\Omega)$ as a subspace of $V$.
Then the part $\mathcal{A}^H$ of $\mathcal{A}$ in $L^2(\Omega)$, i.e., the restriction of $\mathcal{A}$ to
\[
	D(\mathcal{A}^H) \coloneqq \bigl\{ u \in V : \mathcal{A}u \in L^2(\Omega) \bigr\},
\]
acts as an m-accretive operator on $L^2(\Omega)$ by Propositions~\ref{prop:maxmon} and~\ref{prop:maxacc}.
Moreover, the set $D(\mathcal{A}^H)$ is dense in $L^2(\Omega)$ by
Proposition~\ref{prop:dense} since $D(\mathcal{A}^c) \subset D(\mathcal{A}^H)$ and $\mathrm{C}(\overline{\Omega})$
is continuously and densely embedded into $L^2(\Omega)$.
Hence $-\mathcal{A}^H$ generates a non-linear contraction $\mathrm{C}_0$-semigroup $(S(t))_{t \ge 0}$ on $L^2(\Omega)$.
In view of Remark~\ref{rem:generator} we can justly call $u(t,x) \coloneqq (S(t)u_0)(x)$ the
unique solution of~\eqref{eq:ParRobin} for $u_0 \in L^2(\Omega)$,
and we refer to this solution as the $L^2$-solution of~\eqref{eq:ParRobin}. The following is our main result.

\begin{theorem}\label{thm:generation}
	Assume~\eqref{eq:structure2} and~\eqref{eq:monass}. Then $\mathcal{A}^c$ is m-accretive on $\mathrm{C}(\overline{\Omega})$,
	the semigroup generated by $-\mathcal{A}^c$ being the
	restriction of $(S(t))_{t \ge 0}$ to $\mathrm{C}(\overline{\Omega})$.
	Thus for $u_0 \in \mathrm{C}(\overline{\Omega})$ the unique $L^2$-solution $u$ of~\eqref{eq:ParRobin} is continuous
	on the closed parabolic cylinder $[0,\infty) \times \overline{\Omega}$, i.e., continuous up to the parabolic boundary.
\end{theorem}
\begin{proof}
	By Propositions~\ref{prop:maxmon} and~\ref{prop:maxacc} the operator $\mathcal{A}^c$ is m-accretive.
	Hence $-\mathcal{A}^c$ generates
	a non-linear contraction $\mathrm{C}_0$-semigroup $(S^c(t))_{t \ge 0}$ on $\mathrm{C}(\overline{\Omega})$, see Definition~\ref{def:semigroup} and Proposition~\eqref{prop:dense}.
	Since $(I + \alpha \mathcal{A}^c)^{-1}$ is the restriction of $(I + \alpha \mathcal{A}^H)^{-1}$ to $\mathrm{C}(\overline{\Omega})$, see the additional assertion in Proposition~\ref{prop:maxacc},
	the operator $S^c(t)$ is the restriction of $S(t)$ to $\mathrm{C}(\overline{\Omega})$. Thus $t \mapsto S^H(t)u_0$ is continuous
	as a function with values in $\mathrm{C}(\overline{\Omega})$ provided that $u_0 \in \mathrm{C}(\overline{\Omega})$. In this case
	$u(t,x) \coloneqq (S(t)u_0)(x)$ is jointly continuous with respect to $t \ge 0$ and $x \in \overline{\Omega}$.
\end{proof}

\begin{remark}
	Since $\mathcal{A}^c$ is m-accretive, we have a unique solution $u \in \mathrm{C}([0,\infty) \times \overline{\Omega})$
	even for the inhomogeneous problem
	\[
		\left\{ \begin{aligned}
			u_t(t,x) - \divi a(x,\nabla u(t,x)) + b(x,u(t,x)) & = f(t,x) && t > 0, \; x \in \Omega \\
			a(x,\nabla u(t,x)) \cdot \nu + h(x,u(x) & = 0 && t > 0, \; x \in \partial\Omega \\
			u(0,x) & = u_0(x) && x \in \Omega
		\end{aligned} \right.
	\]
	whenever $u_0 \in \mathrm{C}(\overline{\Omega})$ and $f \in L^1_{\mathrm{loc}}( [0,\infty); \mathrm{C}(\overline{\Omega}) )$,
	see~\cite[Corolllary~IV.8.4]{Showalter}.
\end{remark}

Let us finally look at a class of examples that satisfy the conditions~\eqref{eq:structure2} and~\eqref{eq:monass}.
Since in the diffusion equation~\eqref{eq:ParRobin} the function $a(x,\nabla u)$ is the flux, it is natural to assume
that it points into the direction of steepest descent, i.e., into the direction of $\nabla u$.
If we assume in addition that the magnitude of the flux depends only on the steepness of $u$ and possibly
on the location $x \in \Omega$, i.e., if we assume that $a(x,z) = m(x,|z|) \, z$, then we have a
simple criterion to check condition~\eqref{eq:monass}.

\begin{lemma}\label{lem:moncrit}
	Assume that $a(x,z) = m(x,|z|) \, z$ for a measurable function $m\colon \mathds{R} \times \mathds{R}_+ \to \mathds{R}_+$.
	Then the first condition in~\eqref{eq:monass} is satisfied if and only if
	$y \mapsto m(x,y) \, y$ is nondecreasing for every $x \in \Omega$.
	Similar assertions hold for $b$ and $h$.
\end{lemma}
\begin{proof}
	For all $x \in \Omega$ and all $z_1$ and $z_2$ in $\mathds{R}^N$ we have
	\begin{align*}
		& (z_1 - z_2) \, (a(x,z_1) - a(x,z_2)) \\
			& \qquad = m(x,|z_1|) \; |z_1|^2 - \bigl(m(x,|z_1|) + m(x,|z_2|)\bigr) \; z_1 \cdot z_2 + m(x,|z_2|) \; |z_2|^2 \\
			& \qquad \ge m(x,|z_1|) \; |z_1|^2 - \bigl(m(x,|z_1|) + m(x,|z_2|)\bigr) \; |z_1| \; |z_2| + m(x,|z_2|) \; |z_2|^2 \\
			& \qquad = \bigl(|z_1| - |z_2|\bigr) \; \bigl(m(x,|z_1|) \; |z_1| - m(x,|z_2|) \; |z_2|\bigr)
	\end{align*}
	with equality if $z_1$ and $z_2$ point into the same direction.

	If $y \mapsto m(x,y) \, y$ is nondecreasing, then both factors in the last expression have the same sign.
	Thus the product is nonnegative and the first condition in~\eqref{eq:monass} is fulfilled.

	Conversely, if the first condition in~\eqref{eq:monass} is fulfilled, then by choosing $z_2$ to be a positive
	multiple of $z_1$, we see that
	\[
		\bigl(a-b\bigr) \; \bigl(m(x,a) \, a - m(x,b) \, a\bigr) \ge 0
	\]
	for all $a, b > 0$, which implies that $y \mapsto m(x,y) \, y$ is nondecreasing.
\end{proof}

As a consequence of Lemma~\ref{lem:moncrit} we see that an important and commonly encountered class of equations of
$p$-Laplace-type satisfy conditions~\eqref{eq:structure2} and~\eqref{eq:monass}.
For simplicity we consider only examples with $p < N$, but
similar assertions hold true for $p \ge N$ with slightly different integrability assumptions on the coefficients.

\begin{example}
	Let $p \in (1,N)$.
	Set $a(x,z) \coloneqq a_0(x) (s + |z|^{p-2}) \, z$ or $a(x,z) \coloneqq a_0(x) (s^2 + |z|^2)^{\frac{p-2}{2}} \, z$
	with a constant $s \ge 0$ and a measurable function $a_0\colon \Omega \to [\nu,\mu]$, where $0 < \nu \le \mu$.
	Set $b(x,u) \coloneqq b_0(x) |u|^{p-2} u$
	and $h(x,u) \coloneqq h_0(x) |u|^{p-2} u$ with nonnegative measurable functions $b_0 \in L^{\frac{N}{p-\eps}}(\Omega)$
	and $h_0 \in L^{\frac{N-1}{p-1-\eps}}(\partial\Omega)$ for some $\eps > 0$.
	Then the assumption in~\eqref{eq:structure2} and~\eqref{eq:monass} are satisfied and thus Theorem~\ref{thm:generation} applies.
	In the special case $s=0$ and $a_0(x) = 1$ for all $x \in \Omega$ we obtain that the equation
	\[
		\left\{ \begin{aligned}
			u_t(t,x) - \Delta_p u(t,x) + b_0(x) |u|^{p-2} u & = f(t,x) && t > 0, \; x \in \Omega \\
			|\nabla u(t,x)|^{p-2} \tfrac{\partial u(t,x)}{\partial \nu} + h_0(x) |u|^{p-2} u & = 0 && t > 0, \; x \in \partial\Omega \\
			u(0,x) & = u_0(x) && x \in \Omega
		\end{aligned} \right.
	\]
	has a unique $\mathrm{C}^0$-solution $u \in \mathrm{C}( [0,\infty) \times \overline{\Omega} )$.
\end{example}

The strategy of this section applies also to certain dynamic boundary conditions, which are often
called Wentzell-Robin boundary conditions, if we carry out the arguments on a suitable product space,
confer also~\cite{AMPR03} where this idea was originally introduced to the literature for linear equations.
More precisely, we consider the equation
\begin{equation}\label{eq:ParWentzell}
	\left\{ \begin{aligned}
		u_t(t,x) - \divi a(x,\nabla u(t,x)) + b(x,u(t,x)) & = 0 && t > 0, \; x \in \Omega \\
		\beta u_t(t,x) + a(x,\nabla u(t,x)) \cdot \nu + h(x,u(t,x)) & = 0 && t > 0, \; x \in \partial\Omega \\
		u(0,x) & = u_0(x) && x \in \Omega.
	\end{aligned} \right.
\end{equation}
At least formally, we can use the first line in~\eqref{eq:ParWentzell} to express $u_t$ in terms of $u$
in the second line. Then we arrive at what is classically called Wentzell-Robin boundary conditions.

We can show well-posedness of~\eqref{eq:ParWentzell} in the space of continuous functions.
Since the arguments are very similar to what we did before, we only sketch the proof of the following theorem,
but see~\cite{Nittka09} for a detailed proof in the linear case.

\begin{theorem}
	Assume~\eqref{eq:structure2} and~\eqref{eq:monass}. Let $\beta\colon \partial\Omega \to \mathds{R}$ be measurable and such that
	\[
		0 < \essinf_{\partial\Omega} \beta \le \esssup_{\partial\Omega} \beta < \infty.
	\]
	For $u_0 \in \mathrm{C}(\overline{\Omega})$, problem~\eqref{eq:ParWentzell} has a unique solution.
	This solution is continuous on the parabolic cylinder $[0,\infty) \times \overline{\Omega}$.
\end{theorem}
\begin{proof}
	Define the reflexive, separable Banach space
	\[
		\mathcal{V} \coloneqq \bigl\{ (u,u|_{\partial\Omega}) : u \in W^{1,p}(\Omega) \cap L^2(\Omega), \; u|_{\partial\Omega} \in L^2(\partial\Omega) \bigr\}
	\]
	and let $\mathcal{A}_W\colon \mathcal{V} \to \mathcal{V}'$ be defined by the formal expression~\eqref{eq:cAdef}.
	The proof of Proposition~\ref{prop:monop} shows that $\mathcal{A}_W$ is bounded, continuous and monotone.
	For $q \in [2,\infty]$ consider $L^q(\Omega) \oplus_q L^q(\partial\Omega)$, which for $q \in [2,\infty)$ is equipped with the norm
	given by
	\[
		\|(u,g)\|_{L^q(\Omega) \oplus_q L^q(\partial\Omega)}^q \coloneqq \|u\|_{L^q(\Omega)}^q + \|g\|_{L^q(\partial\Omega; \beta \, \dx\sigma)}^q,
	\]
	whereas for $q = \infty$ we set
	\[
		\|(u,g)\|_{L^\infty(\Omega) \oplus_\infty L^\infty(\partial\Omega)} \coloneqq \|u\|_{L^\infty(\Omega)} + \|g\|_{L^\infty(\partial\Omega; \beta \, \dx\sigma)}.
	\]
	Here $\sigma$ denotes the surface measure on $\partial\Omega$, i.e., the $(N-1)$-dimensional Hausdorff measure.
	Set $\mathcal{H} \coloneqq L^2(\Omega) \oplus L^2(\partial\Omega)$, and equip
	$\mathcal{C} \coloneqq \{ (u,u|_{\partial\Omega}) : u \in \mathrm{C}(\overline{\Omega}) \}$
	with the norm of $L^\infty(\Omega) \oplus_\infty L^\infty(\partial\Omega)$. Then both of these spaces
	are subspaces of $V'$ via
	\[
		\apply[V',V]{(u,g)}{(v,v|_{\partial\Omega})} \coloneqq \scalar[L^2(\Omega)]{u}{v} + \scalar[L^2(\partial\Omega;\beta \, \dx\sigma)]{g}{v}
	\]
	for $v \in V$ and $(u,g) \in \mathcal{H}$ or $(u,g) \in \mathcal{C}$, respectively.
	We consider the parts $\mathcal{A}_W^{\mathcal{H}}$ and $\mathcal{A}_W^{\mathcal{C}}$ of $\mathcal{A}_W$ in $\mathcal{H}$
	and $\mathcal{C}$, respectively, with domains
	\begin{align*}
		D(\mathcal{A}_W^{\mathcal{H}}) & \coloneqq \bigl\{ u \in \mathcal{V} : \mathcal{A}_Wu \in \mathcal{H} \bigr\} \\
		D(\mathcal{A}_W^{\mathcal{C}}) & \coloneqq \bigl\{ u \in \mathcal{V} \cap \mathcal{C} : \mathcal{A}_Wu \in \mathcal{C} \bigr\}
	\end{align*}
	Then similar arguments as in the proof of Propositions~\ref{prop:maxmon} and~\ref{prop:maxacc} show that $\mathcal{A}_W^{\mathcal{H}}$
	and $\mathcal{A}_W^{\mathcal{C}}$ are m-accretive on $\mathcal{H}$ and $\mathcal{C}$, respectively. Let $(\mathcal{S}(t))_{t \ge 0}$
	be the semigroup generated by $-\mathcal{A}_W^{\mathcal{H}}$. If $U_0 \in D(\mathcal{A}_W^{\mathcal{H}})$ and $U(t) \coloneqq \mathcal{S}(t)U_0$, then
	$U'(t) + \mathcal{A}_W^{\mathcal{H}} U(t) = 0$ for almost every $t \ge 0$. Writing $U(t) = (u(t),u(t)|_{\partial\Omega})$ this means that
	\begin{align*}
		& \int_\Omega \nabla \eta \; a(x,\nabla u(t)) + \int_\Omega \eta \; b(x,u(t)) + \int_{\partial\Omega} \eta \; h(x,u) \\
			& \quad = \apply{\mathcal{A}_W^{\mathcal{H}}}{(\eta,\eta|_{\partial\Omega})}
			= \scalar[\mathcal{H}]{u'(t)}{(\eta,\eta|_{\partial\Omega})}
			= -\int_\Omega \eta \; u'(t) - \int_{\partial\Omega} \eta \; u'(t) \beta
	\end{align*}
	for all $\eta \in \mathrm{C}^\infty(\overline{\Omega}) \subset \mathcal{V}$ and almost every $t \ge 0$.
	Hence for almost every $t \ge 0$, the function $u(t) \in W^{1,p}(\Omega)$ is a weak solution of~\eqref{eq:ParWentzell}
	with $t$ fixed.
	This justifies that we call the unique $\mathrm{C}^0$-solution of $u'(t) + \mathcal{A}_W^{\mathcal{H}}u(t) = 0$
	with $u(0) = u_0 \in \overline{D(\mathcal{A}_W^{\mathcal{H}})}$, or rather its first component, the (unique) solution of~\eqref{eq:ParWentzell}.

	In the proof of Proposition~\ref{prop:dense} we have seen that
	$\{ u \in W^{1,p}(\Omega) \cap \mathrm{C}(\overline{\Omega}) : \mathcal{A}u \in \mathrm{C}^\infty_c(\Omega) \}$
	is dense in $\mathrm{C}(\overline{\Omega})$, which implies that $\mathcal{A}^{\mathcal{H}}$ and $\mathcal{A}^{\mathcal{C}}$
	are densely defined. Now the same arguments as for Theorem~\ref{thm:generation} show that for every $u_0 \in \mathrm{C}(\overline{\Omega})$
	the unique $L^2$-solution of~\eqref{eq:ParWentzell} is continuous on $[0,\infty) \times \overline{\Omega}$.
\end{proof}

\bibliographystyle{amsplain}
\bibliography{bcnonlinear}

\end{document}